%% file: 2014_sfemVem.tex
\documentclass[times]{nmeauth}
\usepackage{moreverb}

\usepackage{setspace}

\newcommand\BibTeX{{\rmfamily B\kern-.05em \textsc{i\kern-.025em b}\kern-.08em
T\kern-.1667em\lower.7ex\hbox{E}\kern-.125emX}}

\usepackage{url}
\usepackage{graphicx,graphics}
\usepackage{subfigure}
\usepackage{epsfig}
\usepackage{amsmath,amssymb,amsfonts,amsthm}
\usepackage{mathrsfs}
\usepackage{color}
\usepackage{float}
\usepackage{tabularx}
\usepackage{multirow}
\usepackage{stmaryrd}
\usepackage{caption}

\usepackage{pgfplots}
\usetikzlibrary{plotmarks}

\theoremstyle{remark}
\newtheorem{thm}{Theorem}[section]
\newtheorem{rmk}[thm]{Remark}

\newcommand{\Eref}[1]{Equation (\ref{#1})}
\newcommand{\fref}[1]{Figure (\ref{#1})}

\newcommand{\rmd}{\mathrm{d}}
\newcommand{\bveps}{\boldsymbol{\varepsilon}}
\newcommand{\bvsig}{\boldsymbol{\sigma}}

\newcommand{\bigb}{\mathbf{B}}

\newcommand{\dd}{\mathbf{D}}

\newcommand{\kk}{\mathbf{K}}

\newcommand{\cn}{\mathbf{n}}

\newcommand{\qq}{\mathbf{q}}

\newcommand{\uu}{\mathbf{u}}

\newcommand{\xx}{\mathbf{x}}

\def\volumeyear{2013}

\begin{document}

\runningheads{S.~Natarajan, S.~Bordas, E.T.~Ooi}{Equivalence between the SFEM and the VEM}

\title{On the equivalence between the cell-based smoothed finite element method and the virtual element method}

\author{Sundararajan Natarajan$^a$\corrauth, St\'ephane PA Bordas$^{b,c}$, Ean Tat Ooi$^d$}

\address{$^a$Department of Mechanical Engineering, Indian Institute of Technology, Madras, Chennai - 600036, India. \\ $^b$ Facult\'e des Sciences, de la Technologie et de la Communication, University of Luxembourg, Luxembourg. \\$^c$Theoretical and Applied Mechanics, School of Engineering, Cardiff University, Cardiff CF24 3AA, Wales, UK. \\ $^d$School of Science, Information Technology and Engineering, Federation University, Ballarat, VIC 3353, Australia.}

\corraddr{Department of Mechanical Engineering, Indian Institute of Technology, Madras, Chennai - 600036, India. Email: snatarajan@cardiffalumni.org.uk; snatarajan@iitm.ac.in.}

\begin{abstract}
We revisit the cell-based smoothed finite element method (SFEM) for quadrilateral elements and extend it to  arbitrary polygons and polyhedra in 2D and 3D. We highlight the equivalence between the SFEM and the virtual element method (VEM). Based on the VEM, we propose a new stabilization approach to the SFEM when applied to arbitrary polygons and polyhedra. The accuracy and the convergence properties of the SFEM are studied with a few benchmark problems in 2D and 3D linear elasticity. Later, the SFEM is combined with the scaled boundary finite element method to problems involving singularity within the framework of linear elastic fracture mechanics in 2D.
\end{abstract}


\maketitle

\vspace{-20pt}

\section{Introduction}
\label{intro}

The finite element method (FEM) relies on discretizing the domain with non-overlapping regions, called `elements'. In the conventional FEM, the topology of the elements is restricted to triangles and quadrilaterals in 2D or tetrahedrals and hexahedrals in 3D. The use of such standard shapes, simplifies the construction of the approximation over the elements, however, this may require sophisticated (re-) meshing algorithms to either generate high-quality meshes or to capture topological changes. Moreover, the accuracy of the solution depends on the quality of the element employed: Lee and Bathe~\cite{leebathe1993} observed that the shape functions lose their ability to reproduce the displacement fields when the mesh is distorted. In an effort to overcome the limitations of the FEM, research has been focussed on:
\begin{itemize}
\item De-coupling geometry and analysis, for example, meshfree methods~\cite{gingoldmonaghan1977,nguyenrabczuk2008a}, PU enrichment~\cite{melenkbabuvska1996,belytschkogracie2009}, Immersed boundary method~\cite{peskin2002}.
\item Improving the element formulations
\begin{itemize}
\item Strain smoothing~\cite{liunguyen2007,nguyenbordas2008,bordasrabczuk2010,bordasnatarajan2011}
\item Unsymmetric formulations~\cite{rajendranooi2007,rajendran2010}
\item hybrid Trefftz FEM~\cite{szeliu2010,wangqin2011}
\item Polygonal FEM~\cite{sukumarmalsch2006}
\end{itemize}
\item Coupling geometry and analysis, for example, isogeometric analysis~\cite{kaganfischer2003,hughescottrell2006}, isogeometric boundary element method~\cite{scottsimpson2013,simpsonbordas2013}.
\item Boundary based methods, for example, boundary element method~\cite{wrobelaliabadi2002,atroshchenkobordas2014}, scaled boundary finite element method~\cite{wolfsong2000}.
\item Advanced mesh generators~\cite{levyliu2010,yanwang2013,chenwang2014}
\end{itemize}
In this study, we focus and attempt to bridge the gap between two classes of method which both focus on relaxing somewhat the constraints posed on the mesh used in finite element analysis, the strain smoothing technique~\cite{liudai2007} and the virtual element method~\cite{veigamanzini2012}. 

\subsection{Background}

Liu \textit{et al.,}~\cite{liudai2007}, extended the concept of stabilized conforming nodal integration (SCNI)~\cite{chenwang2000} to finite element approximations and coined the resulting method the Smoothed Finite Element Method (SFEM). Liu \textit{et al.} formulated a series of SFEM models: cell-based SFEM (CSFEM)~\cite{liunguyen2007}, node based SFEM (NSFEM)~\cite{liunguyen2009}, edge-based SFEM (ESFEM)~\cite{liunguyen2009a}, face-based SFEM (FSFEM)~\cite{thoiliu2009} and alpha-FEM~\cite{liunguyen2008}. All these SFEM use finite element meshes with linear interpolants, because the idea of the method is to improve the behaviour of simplex elements (triangles, tetrahedra) for which meshes are easier to generate automatically. Only one attempt to use smoothing for higher order approximation is known~\cite{bordasnatarajan2011}, which is the only counter example. In the CSFEM, the elements are divided into smoothing cells, over which the standard (compatible) strain field is smoothed. The method may also be seen as dividing the domain into smoothing domains, which may be constructed independently of the mesh. Yet, from a practical view point, it is simpler to use the mesh data structure to generate the smoothing domain, either from the element interior (cell-based), the nodes (node-based), the edges or face of the elements (edge/face based). This smoothing allows the volume integration to be transformed into a surface integration by employing the divergence theorem and hence the computation of the stiffness matrix requires only information on the boundary of the subcells. It should be noted that only the CSFEM employs quadrilateral elements, whilst all other SFEM models usually rely on simplex elements as reference mesh. When the CSFEM is used with triangular elements, the resulting stiffness matrix is identical to the conventional FEM. The convergence, stability, accuracy and important computational aspects of this method were studied in detail in~\cite{nguyenbordas2008}. The method was further extended to treat various problems in solid mechanics such as plates~\cite{nguyenrabczuk2008}, shells~\cite{nguyen-thanhrabczuk2008}, nearly incompressible elasticity~\cite{ongliu2013,leemihai2014} to name a few. Recently, the strain smoothing method was combined with enrichment methods in~\cite{bordasnatarajan2011,chenrabczuk2012} to model problems with strong discontinuities and singularities. However, certain difficulties still exist as discussed in~\cite{bordasnatarajan2011}.

\begin{figure}[htpb]
\centering
\subfigure[]{\includegraphics[scale=0.35]{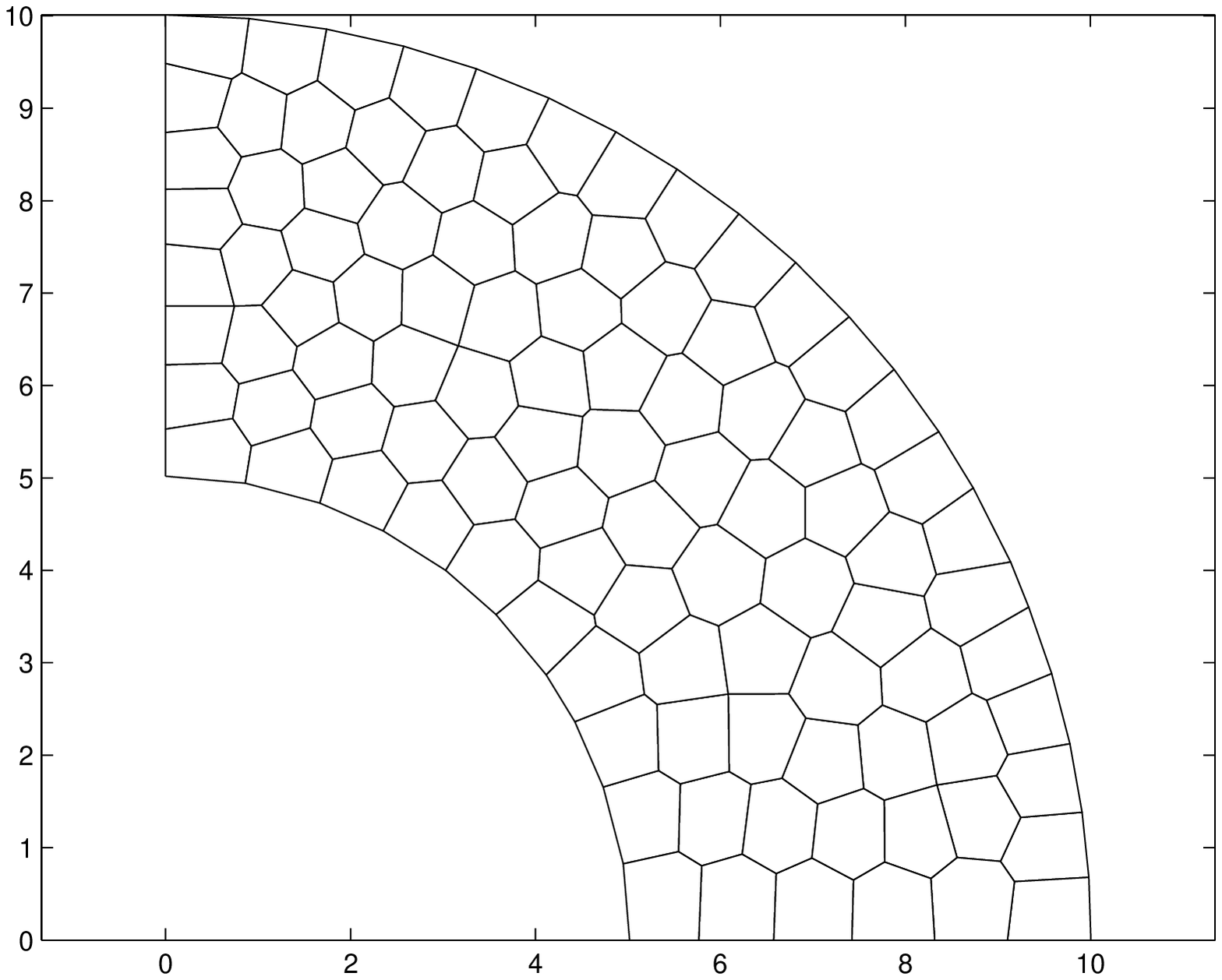}}
\subfigure[]{\includegraphics[scale=0.35]{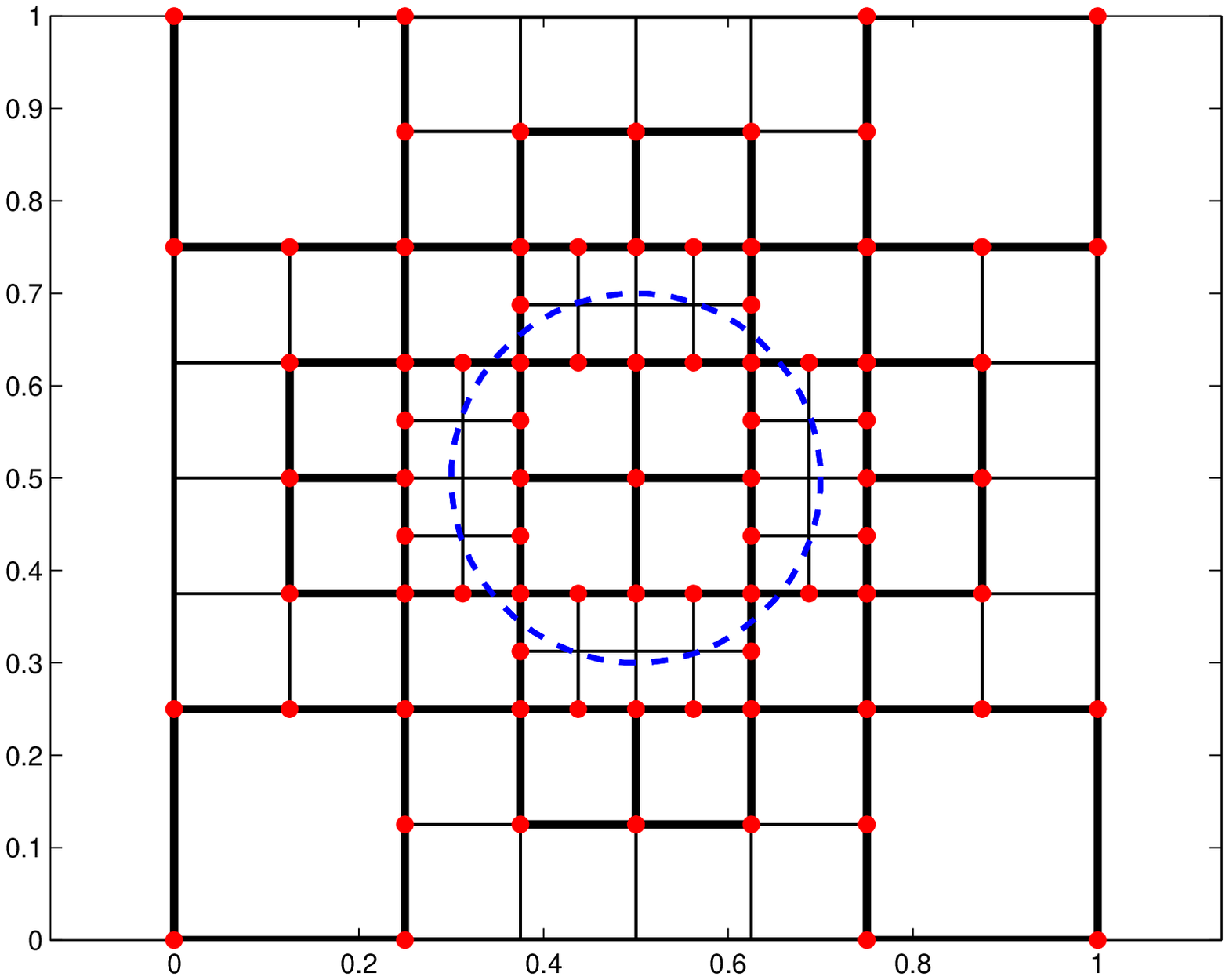}}
\caption{Domain discretized with finite elements: (a) discretization with arbitrary polygons and (b) adaptive refinement leading to a quadtree mesh, where the `dotted' line represents an inner boundary and the `highlighted' elements are the elements with hanging nodes.}
\label{fig:pmeshexample}
\end{figure}

On another related front, polygonal finite element methods (PFEM) have been receiving increasing attention. In PFEM, the domain can be discretized without needing to maintain a particular element topology (see \fref{fig:pmeshexample}). Moreover, this is advantageous in adaptive mesh refinement, where a straightforward subdivision of individual elements usually results in hanging nodes (see \fref{fig:pmeshexample}). Conventionally, this is eliminated by introducing additional edges/faces to retain conformity. This can be avoided if we can directly compute the stiffness matrices on polyhedral meshes with hanging nodes. Polygonal/polyhedral elements allow to treat all elements within a quad-tree/oct-tree mesh within a single paradigm. For example, elements of class (1) quadrilaterals, (2) pentagons and (3) hexagons can be assembled within a single routine.

In 1971, Wachspress~\cite{wachspress1971} developed a method based on rational basis approximations for elements with an arbitrary number of sides. However, these elements were not easily used because of difficulties in constructing the basis functions, performing numerical integration, but more importantly generating polyhedral meshes on arbitrary domains.. Thanks to advancements in mathematical software, viz., Mathematica \textsuperscript{\textregistered} and Maple \textsuperscript{\textregistered} and the pioneering work of Alwood and Cornes [15], Sukumar and Tabarraei [16], Dasgupta [17], to name a few and the recent drive from the computer graphics/science community~\cite{levyliu2010,yanwang2013,chenwang2014}, the discretization of the domain with finite elements having arbitrary number of sides/edges is possible for realistic application. 

Once the mesh is generated and the basis functions are constructed, the conventional Galerkin procedure is normally employed to solve the governing equations over the polygonal/polyhedral meshes. However, the numerical integration of the terms in the stiffness matrix over polygonal/polyhedral meshes requires special techniques and often is complicated. Improving numerical integration over polytopes has thus beers the focus of significant attention in the recent literature~\cite{sukumartabarraei2004,natarajanbordas2009,mousavixiao2010,talischipaulino2013,sudhakaralmeida2014}. The strain smoothing technique is another alternative to integrate over arbitrary polygons. Dai \textit{et al.,}~\cite{dailiu2007} observed that on an arbitrary polygon with $n >$ 4 (where $n$ is the number of sides of the polygon) must have a minimum of $n$ subcells to ensure stability. Recently, Natarajan~\textit{et al.,}~\cite{natarajanooi2014} and Ooi \textit{et al.,}~\cite{ooisong2012} employed the scaled boundary finite element method over arbitrary polygons for linear elasticity and linear elastic fracture mechanics, respectively. 

The virtual element method (VEM) was recently discovered in~\cite{veigamanzini2012,daveigabrezzi2013,daveigabrezzi2014,manzinirusso2014}. The VEM has evolved from the mimetic difference methods~\cite{veigamanzini2012}. It is a Galerkin approach, however, unlike the conventional FEM, the VEM does not require an explicit form for the basis functions to compute the stiffness matrix. Moreover, the VEM can be employed over arbitrary polygons and polyhedra. Within the framework of the VEM, the space within an element of decomposition\footnote{There is no restriction on the shape of the element in VEM}, contains certain polynomials that guarantee \textit{accuracy} and additional functions for stability. The VEM also alleviates the numerical integration difficulty encountered in the conventional polygonal FEM. As the method does not require the knowledge of the shape functions in explicit form, the implementation is computationally less intensive.

\subsection{Objective}
We saw above that the strain smoothing technique can be applied to arbitrary polygons. However, it was observed in~\cite{natarajanooi2014} that the strain smoothing technique over arbitrary polytopes yield less accurate solution when compared to other techniques, such as the conventional polygonal finite element method. Moreover, like the conventional finite elements, the polygonal elements requires special techniques, such as enrichment methods to treat problems with strong discontinuities and singularities. On the other hand, the virtual element method enables accurate and stable solution for completely arbitrary polygons. It therefore seems sensible to investigate the eventual connections between both classes of methods. More specifically, the main objectives of the paper are:

\begin{itemize}
\item to revisit the strain smoothing technique, in particular, the cell based SFEM and extend it to arbitrary polygons and polyhedra. 
\item investigate connections between the strain smoothing technique and the VEM.
\item propose a new stabilization for the SFEM with one subcell inspired from the VEM, and which is applicable to arbitrary polygons.
\item to study the accuracy and the convergence properties of the SFEM with the new stabilization technique.
\item to couple the SFEM with the scaled boundary finite element method (SBFEM) (see Section \ref{sfemsbfem} for a detailed discussion on the SBFEM) to study problems with strong discontinuities and singularities.
\end{itemize}
Throughout this paper, SFEM stands for the cell-based smoothed finite element method unless mentioned otherwise.

\subsection{Outline} The paper is organized as follows. Section \ref{csfemintroduction} revisits the basics of the cell-based smoothed finite element method as applied to 2D and 3D elasticity. Section \ref{vem} briefly reviews the virtual element method. The similarity between the SFEM and the VEM is discussed in Section \ref{sfemvsvem}. Some numerical examples are analyzed to demonstrate the accuracy and the convergence properties of the cell-based smoothed finite element method in Section 5 with a few problems taken from linear elasticity. The SFEM is combined with the SBFEM and the accuracy of the approach is demonstrated through benchmark problems in linear elastic fracture mechanics. The major conclusions and future research directions are summarised in the final section.

\section{Overview of the smoothed finite element method}
\label{csfemintroduction}
\input{sfemintro}

\section{Basics of the virtual element method}
\label{vem}
\input{vemintro}

\section{Equivalence between the SFEM and the VEM}
\label{sfemvsvem}
In this section, we shall attempt to demonstrate the equivalence between the SFEM and the VEM, first by presenting a generic expression and then by computing the stiffness matrix of a generic polygonal element. We do this for the following model problem in two dimensions:
\begin{equation}
\Delta u = f \hspace{3pt} \textup{in} \hspace{2pt} \Omega; \hspace{4pt} u = 0 \hspace{3pt} \textup{on} \hspace{5pt} \partial \Omega
\end{equation}
where $\Omega \in \mathbb{R}^2$. The corresponding variational formulation reads: $\textup{find}~~ u \in \mathcal{V} \colon= H_0^1(\Omega)$ such that:
\begin{equation}
\forall v \in \mathcal{V} \hspace{10pt} a(u,v) = (f,v) 
\end{equation}
where $(\cdot,\cdot)$ represents the scalar inner product in $L^2$ and $a(u,v) = (\nabla u, \nabla v)$. The VEM for the above bilinear form starts by defining a projection operator that follows the orthogonality condition~\cite{daveigabrezzi2014}:
\begin{equation}
\forall p_k \in \mathcal{P}(E) \hspace{10pt} \left( \nabla p_k, \nabla( \Pi^\nabla v_h - v_h) \right) = 0
\label{eqn:orthocondoperator}
\end{equation}
where $v_h \in V(E)$ is the finite element space and $\mathcal{P}(E)$ is the space of polynomials over a polygonal element $E$ which has the basis $m_\alpha$, where
\begin{equation*}
m_\alpha \colon= \left( \frac{\mathbf{x}-\mathbf{x}_\mathcal{D}}{h_\mathcal{D}} \right)^\alpha
\end{equation*}
where $\boldsymbol{\alpha}  = (\alpha_1,\alpha_2), |\boldsymbol{\alpha}| \le k$, $k$ is the degree of the polynomial and $\mathbf{x} = (x,y)$, $\mathbf{x}_\mathcal{D}$ and $h_\mathcal{D}$ is the centroid and the diameter of the element, respectively. When $k=$ 1, \Eref{eqn:orthocondoperator}, becomes~\cite{daveigabrezzi2014}
\begin{align}
\nabla p_1 \cdot \nabla(\Pi^\nabla v_h) &= \frac{1}{|E|} \nabla p_1 \cdot \int\limits_E \nabla v_h \nonumber \\
&= \frac{1}{|E|} \int\limits_E \nabla v_h = \nabla(\Pi^\nabla v_h) \colon= \mathbf{g}(v_h)
\end{align}
where we have used $p_1 = x_1$ and $p_1=x_2$. Hence,
\begin{equation}
\Pi^\nabla v_h = \mathbf{x} \cdot \mathbf{g}(v_h) + \textup{constant}
\label{eqn:operatorterm}
\end{equation}
The local stiffness matrix of the virtual element method is written as:
\begin{equation*}
(\kk^{\rm VEM})_{IJ} = (\nabla \Pi^\nabla \phi_I,\nabla \Pi^\nabla \phi_J) + \textup{additional terms for stability}
\end{equation*}
The first term in the above equation is called the \emph{consistency} term, which by using \Eref{eqn:operatorterm} can be written as:
 \begin{equation*}
 (\nabla \Pi^\nabla \phi_I,\nabla \Pi^\nabla \phi_J)  = |E| ~\mathbf{g}(\phi_I) \cdot \mathbf{g}(\phi_J)
 \end{equation*}
 when $\phi_I$ is linear on the edge of the element, 
 \begin{equation}
 \mathbf{g}(\phi_I) = \frac{1}{|E|} \int\limits_E \nabla \phi_I = \frac{1}{2|E|} \left( \ell_{I-1} \mathbf{n}_{I-1} + \ell_{I} \mathbf{n}_{I} \right)
 \label{eqn:vemgradientterm}
 \end{equation}

In case of the SFEM, the smoothed gradient matrix is given by:
\begin{align}
\bigb &= \left[ \begin{array}{cc} \frac{\partial N_I}{\partial x} \frac{\partial N_I}{\partial y}  \end{array} \right] \nonumber \\
\tilde{\bigb} &= \frac{1}{A_C} \int\limits_{S_C} \mathbf{n}^{\rm T} N_I(\xx)~\rmd S
\end{align}
where the integration is over the subcell boundary. When $S_C=$ 1 and when the shape functions are linear on the boundary, the gradient matrix is given by:
\begin{align}
\renewcommand{\arraystretch}{3}
\tilde{\bigb} &= \frac{1}{2A_C} \left\{ \begin{array}{c} \ell_{I-1} n^x_{I-1} + \ell_{i} n^x_{I} \\\ell_{I-1} n^y_{I-1} + \ell_{I} n^y_{I} \end{array} \right\} \nonumber \\
&= \frac{1}{2A_C} \left( \ell_{I-1} \mathbf{n}_{I-1} + \ell_{I} \mathbf{n}_{I} \right)
\end{align}
By comparing the above equation with \Eref{eqn:vemgradientterm}, we can easily see that the gradient matrix obtained by the smoothing technique is identical to the \emph{consistency} term of the VEM. 

Next, we demonstrate the equivalence between the methods with two worked examples by computing the stiffness matrix of (a) a quadrilateral element and (b) a pentagon. For comparison, we also present the stiffness matrix computed by conventional FEM by employing Gaussian quadrature. The following expressions are used for the respective methods, viz., the FEM, the SFEM and the VEM, to compute the stiffness matrix:
\begin{itemize}
\item $\kk^{\rm FEM} = \int\limits_\Omega \bigb^{\rm T} \bigb~\rmd \Omega$.
\item $\kk^{\rm SFEM} = \sum\limits_{C=1}^{nc} \int\limits_{\Omega_C} \tilde{\bigb}^{\rm T}_C \tilde{\bigb}_C~\rmd \Omega$
\item $\kk^{\rm VEM} = \frac{ \mathbf{R} \mathbf{R}^{\rm T}}{|E|} + ( \mathbf{I}-\mathbf{\Pi})^{\rm T} ( \mathbf{I}-\mathbf{\Pi})$
\end{itemize}
where $\bigb$ is the gradient matrix and $\tilde{\bigb}$ is the smoothed gradient matrix given by:
\begin{align}
\bigb &= \left[ \begin{array}{cc} \frac{\partial N_I}{\partial x} \frac{\partial N_I}{\partial y}  \end{array} \right] \nonumber \\
\tilde{\bigb} &= \frac{1}{A_C} \int\limits_{S_C} \mathbf{n}^{\rm T} N_I(\xx)~\rmd S
\end{align}
and $\mathbf{\Pi} = \tilde{\mathbf{\Pi}} + \mathbf{\Pi}_o ( \mathbf{I} - \tilde{\mathbf{\Pi}}  )$ and $ |E| \tilde{\mathbf{\Pi}}  = \mathbf{N} \mathbf{R}^{\rm T}$. The matrices $\mathbf{R}$ and $\mathbf{N}$ are given by~\cite{manzinirusso2014,daveigabrezzi2014}:
\begin{equation*}
\mathbf{R} = \frac{1}{2} \left[ \begin{array}{c} \ell_n \mathbf{n}_n + \ell_1 \mathbf{n}_1 \\ \ell_1 \mathbf{n}_1 + \ell_2 \mathbf{n}_2 \\ \cdots \\ \ell_{n-1} \mathbf{n}_{n-1} + \ell_n \mathbf{n}_n \end{array} \right] \hspace{10pt} \mathbf{N} = \left[ \begin{array}{cc} x_1 & y_1 \\ x_2 & y_2 \\ \cdots \\ x_n & y_n \end{array} \right]
\end{equation*}
where $\ell_i (i=1,\cdots,n)$ is the length of edge $i$, $x_n,y_n$ are the coordinates of the vertex of the polygon, $|E|$ is the measure of the polygon.

\subsection{Stiffness matrix for the unit square} In this case, consider a unit square $[0,1]\times[0,1]$. We compute the stiffness matrix using the conventional FEM with bilinear shape functions, the VEM and the SFEM with both one and several subcells.

\paragraph{Finite element} The stiffness matrix computed from the classical bilinear finite elements with reduced integration (i.e., one Gauss point at the center of the element) and full integral (four Gauss points):
\begin{equation}
\kk^{\rm FEM}_ {\rm red} = \frac{1}{2} \left[ \begin{array}{rrrr} 1 & 0 & -1 & 0 \\ 0 & 1 & 0 & -1 \\ -1 & 0 & 1 & 0 \\ 0 & -1 & 0 & 1 \end{array} \right] \hspace{15pt}
\kk^{\rm FEM}_ {\rm full} = \frac{1}{12} \left[ \begin{array}{rrrr} 8 & -2 & -4 & -2 \\ -2 & 8 & -2 & -4 \\ -4 & -2 & 8 & -2 \\ -2 & -4 & -2 & 8 \end{array} \right] 
\end{equation}

\paragraph{Virtual element} In this case, we use the order of the monomial $k=$ 1. This implies that the shape functions on the boundary of the elements are linear. As noted in Section \ref{vem}, the stiffness matrix computed from the VEM has two parts (see \Eref{eqn:vem3d}): (a) the consistency term and (b) the stability term. The consistency and the stability terms are given by:
\begin{equation}
\kk^{\rm VEM}_ {\rm const} = \frac{1}{2} \left[ \begin{array}{rrrr} 1 & 0 & -1 & 0 \\ 0 & 1 & 0 & -1 \\ -1 & 0 & 1 & 0 \\ 0 & -1 & 0 & 1 \end{array} \right] \hspace{15pt}
\kk^{\rm VEM}_ {\rm stab} = \frac{1}{4} \left[ \begin{array}{rrrr} 1 & -1 & 1 & -1 \\ -1 & 1 & -1 & 1 \\ 1 & -1 & 1 & -1 \\ -1 & 1 & -1 & 1 \end{array} \right] 
\end{equation}
and the final stiffness matrix is computed by adding the consistency term and the stability term:
\begin{equation}
\kk^{\rm VEM} = \frac{1}{12} \left[ \begin{array}{rrrr} 9 & -3 &  -3 & -3 \\ -3 & 9 & -3 & -3 \\ -3 & 9 & -3 & -3 \\ -3 & -3 & -3 & 9 \end{array} \right] 
\end{equation}

\paragraph{SFEM} In this case, we use one and two subcells. The stiffness matrices with one subcell and two subcells are given by:
\begin{equation}
\kk^{\rm SFEM}_{\rm SC1Q4}  = \frac{1}{2} \left[ \begin{array}{rrrr} 1 & 0 & -1 & 0 \\ 0 & 1 & 0 & -1 \\ -1 & 0 & 1 & 0 \\ 0 & -1 & 0 & 1 \end{array} \right] \hspace{15pt}
\kk^{\rm SFEM}_{\rm SC2Q4} = \frac{1}{16} \left[ \begin{array}{rrrr} 9 & -1 & -7 & -1 \\ -1 & 9 & -1 & -7 \\ -7 & -1 & 9 & -1 \\ -1 & -7 & -1 & 9 \end{array} \right]
\end{equation}

\subsection{Stiffness matrix for the pentagon} The coordinates of the pentagon are: $[ (0,0), (3,0), (3,2), (3/2,4), (0,4)]$. The consistency and the stability term for the VEM are given by:
\begin{align}
\kk^{\rm VEM}_{\rm const} &= \left[ \begin{array}{rrrrr}  0.5952 &  0.0238 &   -0.4881  &  -0.4048 &    0.2738 \\
0.0238  &   0.3095 &    0.0833&    -0.1190 &   -0.2976 \\
   -0.4881 &    0.0833&     0.4345 &    0.2976&    -0.3274 \\
   -0.4048 &   -0.1190  &   0.2976 &    0.3095  &  -0.0833 \\
    0.2738&    -0.2976  &  -0.3274 &   -0.0833 &    0.4345 \end{array} \right] \nonumber \\
\kk^{\rm VEM}_{\rm stab} &= \left[ \begin{array}{rrrrr} 0.7422  & -0.1966&   -0.3412 &  -0.2578 &  0.0534 \\
   -0.1966 &   0.7422&   -0.3412  & -0.1354 &  -0.0690 \\
   -0.3412 &  -0.3412 &   0.9896 &   0.0364  & -0.3437 \\
   -0.2578&   -0.1354 &   0.0364&    0.8646 &  -0.5078 \\
    0.0534 &  -0.0690&   -0.3437 &  -0.5078&    0.8672 \end{array} \right]
\end{align}
and the stiffness matrix computed by employing the smoothing technique over the pentagon with one subcell is:
\begin{align}
\kk^{\rm SFEM}_{\rm one cell} &= \left[ \begin{array}{rrrrr}  0.5952 &  0.0238 &   -0.4881  &  -0.4048 &    0.2738 \\
0.0238  &   0.3095 &    0.0833&    -0.1190 &   -0.2976 \\
   -0.4881 &    0.0833&     0.4345 &    0.2976&    -0.3274 \\
   -0.4048 &   -0.1190  &   0.2976 &    0.3095  &  -0.0833 \\
    0.2738&    -0.2976  &  -0.3274 &   -0.0833 &    0.4345 \end{array} \right]
\end{align}

From the two examples, presented above, it can be observed that the consistency term of the VEM, stiffness matrix using SFEM with one subcell coincides with the conventional FEM with reduced integration. However, after the addition of the stability term in the VEM or increasing the number of subcells in the SFEM, we observe that the stiffness matrix computed from these approaches are different. It is also noted that the consistency term of the VEM is similar to the SFEM with one subcell. This observation is also true in the case of the pentagon.

The SFEM starts with an assumption that the strain is constant within the subcell and then employs the divergence theorem to convert the domain integral into a surface integral. This suppresses the need to compute the derivatives of the shape functions and the stiffness matrix is computed from the information available on the boundary. When linear elements are employed on the boundary, this assumption holds true. However, for higher order elements, this assumption breaks down. This was observed in~\cite{bordasnatarajan2011}, when the strain smoothing technique was used for Q8 and Q9 and enriched approximations.

In the case of the VEM, no such assumption of constant strain is made over the element. However, the method starts by assuming the variation of the shape functions on the boundary of the element. The method then employs the divergence theorem after defining the projection operators~\cite{veigabrezzi2013,daveigabrezzi2014}. Hence, when linear variation is assumed, the VEM with consistency term and the SFEM with one subcell coincide. This is also true when the number of sides is greater than 4, as noted in the previous example. The stiffness matrix of the VEM has two parts: (a) the first term ensures \textit{consistency} and this term must be computed exactly and (b) the second term ensures \textit{stability}, this can be approximated. The important features of the stability term are: (a) it should scale like the consistency term and (b) should be positive definite. Different choices of stability terms are possible as discussed in~\cite{veigabrezzi2013,veigabrezzi2013a,gaintalischi2014}.  However, for this study, we employ the following stability term based on the work of Beir\~ao Da Veiga \textit{et al.,}~\cite{veigabrezzi2013,veigabrezzi2013a}:
\begin{equation}
\kk_2 = \alpha \mathbf{P}
\label{eqn:stabterm}
\end{equation}
where $\alpha = \alpha^\ast \textup{trace}(\kk^{\rm VEM}_{\rm const})$ and $\mathbf{P}$ is the orthogonal projection operation and is chosen as:
\begin{equation}
\mathbf{P} = \mathbf{I} - \mathbf{T} \left( \mathbf{T}^{\rm T} \mathbf{T} \right)^{-1} \mathbf{T}^{\rm T}
\end{equation}
and the matrix $\mathbf{T}$ is the modified nodal coordinate matrix $N$ of dimension 3$n\times$ 12 in the case of 3D:
\begin{equation}
T_{3I-2 \colon 3I} = \left[ \begin{array}{rrrrrrrrrrrr} 1 & 0 & 0 & y_I & 0 & -z_I & x_I & 0 & 0 & y_I & 0 & z_I\\ 0 & 1 & 0 & -x_I & z_I & 0 & 0 & y_I & 0 & x_I & z_I & 0 \\ 0 & 0 & 1 & 0 & -y_I & x_I & 0 & 0 & z_I & 0 & y_I & x_I \end{array} \right]
\end{equation}
where $\alpha^\ast$ is a scaling coefficient chosen based on a parametric study conducted in the next section. We conclude that the CSFEM is a special case of the more general VEM. Instead of increasing the number of subcells, we add to the one-subcell, the stability term borrowed from the VEM. In this present study, we employ the following form for the stiffness matrix:
\begin{equation}
\kk^{h} = \kk_1 + \kk_2
\end{equation}
where
\begin{equation}
\kk_1 = \int\limits_{\Omega_C} \tilde{\bigb}^{\rm T}_C \dd \tilde{\bigb}_C~\rmd \Omega
\end{equation}
is computed by employing the strain smoothing technique and
\begin{equation}
\kk_2 = \alpha \mathbf{P}
\label{eqn:stabilitytermvem}
\end{equation}
Since the consistency term of the VEM and the SFEM with one subcell are identical, we have $\alpha = \alpha^\ast \textup{trace}(\kk_1)$. This approach enables the use of the most attractive of all cell-based smoothed finite element method (with one subcell), whilst, guaranteeing stability.

\section{Numerical examples}
\label{numres}
In the first part of this section, we employ the cell-based smoothed finite element method with the new stabilization approach to two-dimensional benchmark problems in linear elasticity. The results from the new approach are compared with analytical solution where available and with the conventional FEM. The SFEM with stabilization is applied to polygonal elements in 2D and the accuracy and the convergence properties are studied in detail. Later, the proposed SFEM with stabilization is extended to 3D problems with hexahedral and polyhedral elements. Again, the accuracy and the convergence properties of the proposed method are studied with a patch test and a cantilever beam loaded in shear. In the last part of the section, the SFEM is combined with the scaled boundary FEM for problems involving strong discontinuity and singularities. The results are compared with available solutions in the literature. 

The built-in Matlab \textsuperscript{\textregistered} function {\small voronoin} and Matlab \textsuperscript{\textregistered} functions in PolyTop~\cite{talischipaulino2012} for building mesh-connectivity are used to create the polygonal meshes. For polyhedra meshes, the open source software Neper~\cite{queydawson2011} is employed for building the mesh-connectivity. For the purpose of error estimation and convergence studies, the  error, $L^2$ and $H^1$ norms are used. The displacement norm is given by:
\begin{equation}
|| \uu - \uu^h ||_{L^2(\Omega)} = \sqrt{ \int\limits_\Omega \left[ (\uu-\uu^h) \cdot (\uu - \uu^h) \right]~\rmd \Omega}
\end{equation}
where $\uu^h$ is the numerical solution and $\uu$ is the analytical or a reference solution. The energy norm is given by:
\begin{equation}
|| \uu - \uu^h ||_{H^1(\Omega)} = \sqrt{ \int\limits_\Omega \left[ (\bveps-\bveps^h) \dd (\bveps - \bveps^h) \right]~\rmd \Omega}
\end{equation}

\subsection{Applications to two dimensional problems}
\input{twoDexamples}

\subsection{Applications in three dimensional problems}
\input{threeDexamples}

\input{ApplicationLefm}

\section{Conclusions}
In this paper, we revisited the cell-based smoothed finite element method (SFEM) and constructed a one subcell polygonal/polyhedral smoothed finite element which is stable, accurate and convergent compared to existing approaches. We also demonstrated the equivalence of the method with the virtual element method (VEM). We conclude that the SFEM can be seen as a \textit{special case} of the more general VEM. By utilizing the concept of the stabilizing term from the VEM, we proposed a new stabilized smoothed finite element method. Instead of increasing the number of subcells, we add to the one-subcell, the stability term borrowed from the VEM. When applied to arbitrary polygons/polyhedra, with the proposed method, sub-triangulation of the polygon/polyhedra is not required to ensure stability. From the detailed numerical study, we observe that the SFEM with stabilization term yields more accurate results than the conventional SFEM with many subcells for much fewer numerical operations. To study problems with singularities, the proposed was combined with the scaled boundary finite element method. The proposed method is flexible, easy to implement and yields accurate results.

\ack The authors would also like to acknowledge Dr. Romain Quey for his generous support with the Neper software. Stéphane Bordas would like to thank the support from the European Research Council Starting Independent Research Grant (ERC Stg grant agreement No. 279578) entitled “Towards real time multiscale simulation of cutting in non-linear materials with applications to surgical simulation and computer guided surgery as well as partial support from the EPSRC under grant EP/G042705/1 Increased Reliability for Industrially Relevant Automatic Crack Growth Simulation with the eXtended Finite Element Method and EP/I006494/1 Sustainable domain-specific software generation tools for extremely parallel particle-based simulations.

Simulations were supported by ARCCA and High Performance Computing (HPC) Wales, a company formed between the Universities and the private sector in Wales which provides the UK’s largest distributed supercomputing network.

\bibliographystyle{wileyj}
\bibliography{vem}

\end{document}

%% file: sfemintro.tex
\subsection{Background}

\begin{figure}[htpb]
\centering
\scalebox{0.5}{\input{./Figures/sfemid.pstex_t}}
\caption{Subdivision of an element into subcells: (a) quadrilateral element and (b) arbitrary polygon.}
\label{fig:subcellRepre}
\end{figure}
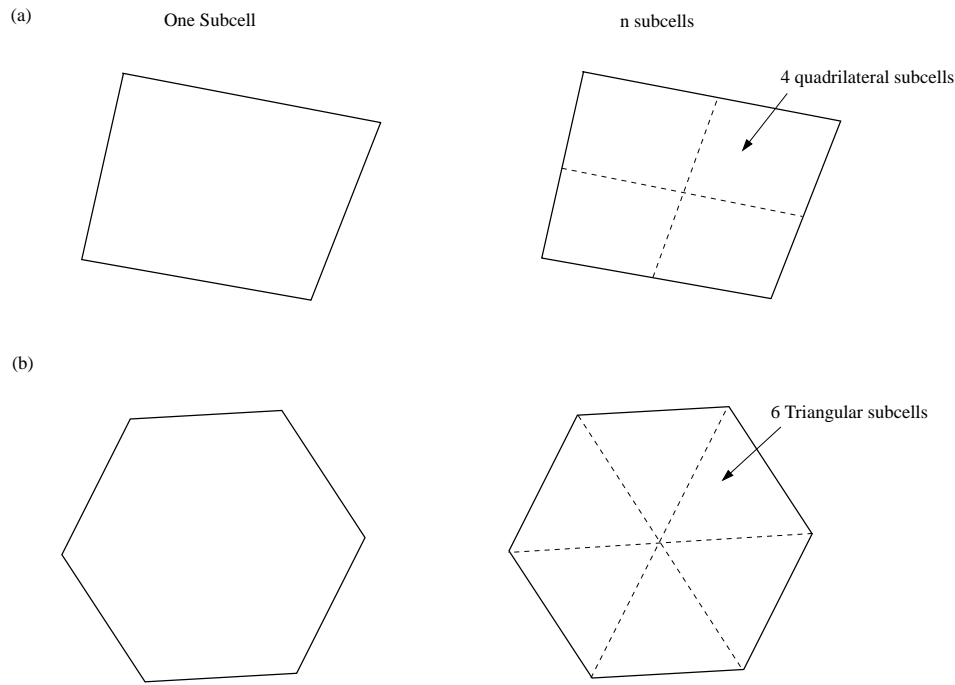

The strain-smoothing method (SSM) was proposed in~\cite{chenwang2000} where the strain is written as the divergence of a spatial average of the standard (compatible) strain field --i.e. symmetric gradient of the displacement field. In the cell-based SFEM, the elements are divided into subcells as shown in \fref{fig:subcellRepre}. The strain field $\tilde {\varepsilon }_{ij}^h$, used to compute the stiffness matrix is computed by a weighted average of the standard strain field $ {\varepsilon }_{ij}^h$. At a point $\xx_C$ in an element $\Omega^h$,

\begin{figure}[htbp]
\centering
\includegraphics[angle=0,width=0.5\textwidth]{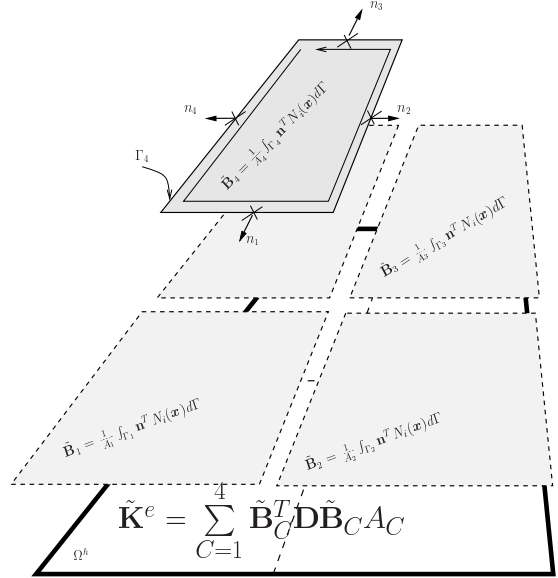}
\caption{Calculation of the smoothed discretized gradient operator.}
\label{fig:bsmoothed}
\end{figure}
\begin{equation}
\tilde {\varepsilon }_{ij}^h (\xx_{C} )=\int_{\Omega ^h}
{\varepsilon _{ij}^h (\xx)\Phi (\xx-\xx_{C} ) \rmd \xx }
\label{eqn:epsilonvar}
\end{equation}
where $\Phi $ is a smoothing function that generally satisfies the following properties~\cite{yoomoran2004} 
\begin{equation}
\Phi \geq 0 ~~~~ \text{and} ~~~~ \int_{\Omega ^h} {\Phi(\xx) \rmd \xx
}=1 \label{eqn:propersmoothfunc}
\end{equation}
\begin{equation}
\Phi = \frac{1}{A_C} ~~~~ \text{in} ~~~~ \Omega_C \quad \text{and}
\quad \Phi = 0 \quad \text{elsewhere}
 \label{eqn:propersmoothfunc2}
\end{equation}
To use~\Eref{eqn:epsilonvar}, the subcell containing point $\xx_C$ must first be located in order to compute the correct value of the weight function $\Phi$. 

The discretised strain field is computed, through the so-called smoothed discretised gradient operator $\tilde{\bigb}$, defined by (see~\fref{fig:bsmoothed} for a schematic representation of the construction)
\begin{equation}
\tilde{\varepsilon}^h(\xx_C)=\tilde{\bigb}_C(\xx_C)\qq
\end{equation}
where $\qq$ contains unknown displacements coefficients defined at a node of a finite element. This definition is similar to the conventional FEM. The smoothed element stiffness matrix for element $e$ is computed by the \emph{sum of the contributions of the subcells} (\fref{fig:bsmoothed})\footnote{The subcells $\Omega_C$ form a partition of the element $\Omega^h$.}
\begin{equation}
\tilde {\kk}^e =\sum\limits_{C=1}^{nc}\int_{\Omega_C} {\tilde
{\bigb}_C^{\textup T} \dd\tilde {\bigb}_C } \rmd \Omega =\sum\limits_{C=1}^{nc}
{\tilde {\bigb}_C^T \dd\tilde {\bigb}_C } \int_{\Omega_C} \rmd \Omega
=\sum\limits_{C=1}^{nc}{\tilde {\bigb}^T_C \dd\tilde {\bigb}_CA_C }
\label{eqn:stiffnessvar}
\end{equation}
where $nc$ is the number of the smoothing cells of the element. The strain displacement matrix $\tilde \bigb_{C}$ is constant over each $\Omega_C$ and is of the following form
\begin{equation}
\tilde \bigb_{C}=\left[ {\begin{array}{*{20}c} \tilde {\bigb}_{C1} &
\tilde {\bigb}_{C2} & \tilde {\bigb}_{C3} & \cdots \tilde {\bigb}_{Cn}
\end{array}} \right] \label{eq:Btilde3dbis}
\end{equation}
where for all shape functions $I \in \{1,\dots,n\}$, the $3 \times 2$ submatrix $\tilde {\bigb}_{CI}$ represents the contribution to the strain displacement matrix associated with shape function $I$ and cell $C$ and  writes  (see~\fref{fig:bsmoothed})
\begin{align}
\renewcommand\arraystretch{2}
\forall I\in\{1,2,\dots,n\},\forall C \in \{1,2,\dots nc \}\tilde
{\bigb}_{CI} &=\frac{1}{A_C} \int_{S_C}\cn^T(\xx) N_I(\xx) \rmd S \nonumber \\ &=
\frac{1}{A_C} \int_{S_C} \left[
\begin{array}{*{20}c}
 { n_x} & 0 \\
 0 & { n_y} \\
 { n_y} & { n_x} \\
\end{array}  \right](\xx) N_I(\xx) \rmd S
\label{eqn:Btilde2d}
\end{align}
Note that since~\Eref{eqn:Btilde2d} is computed on the boundary of $\Omega_C$ and one Gau\ss{} point is sufficient for an exact integration:
\begin{equation}
\renewcommand\arraystretch{2}
\tilde {\bigb}_{CI}(\xx_C)=\frac{1}{A_C}\sum\limits_{b=1}^{nb} {\left(
{{\begin{array}{*{20}c}
 {N_I\left(\xx_b^G\right) n_x} & 0 \\
 0 & {N_I\left(\xx_b^G\right) n_y} \\
 {N_I\left(\xx_b^G\right) n_y} & {N_I\left(\xx_b^G\right) n_x} \\
\end{array} }} \right)} l_b^C
\label{eqn:equation3.14}
\end{equation}
where $\xx_b^G$ and ${l}_b^C$ are the center point (Gau\ss{} point) and the length of $\Gamma_b^C$, respectively. Until now, no assumption was made on the shape of the element. The procedure outlined so far is general and is applicable to polygons of arbitrary shapes~\cite{dailiu2007,natarajanooi2014}. Due to the process of strain smoothing, only the shape function is involved in the calculation of the field gradients and hence the stiffness matrix. In this study, we employ the simple averaging technique to compute the shape functions over arbitrary polygons. The construction of shape function is as follows: for a general polygonal element, the central point $O$ is located by:
\begin{equation}
(x_o,y_o) = \frac{1}{n} \sum\limits_i^n(x_i,y_i)
\end{equation}
where $n$ is the number of nodes of the polygonal element. The shape function at point $O$ is given by $[1/n~\cdots~1/n]$ with size 1 $\times n$. 

\subsection{Extension to 3D}
\begin{figure}[htpb]
\centering
\scalebox{0.5}{\input{./Figures/sfemid3d.pstex_t}}
\caption{Subdivision of an element into subcells: (a) quadrilateral element and (b) arbitrary polygon.}
\label{fig:subcellRepre3d}
\end{figure}
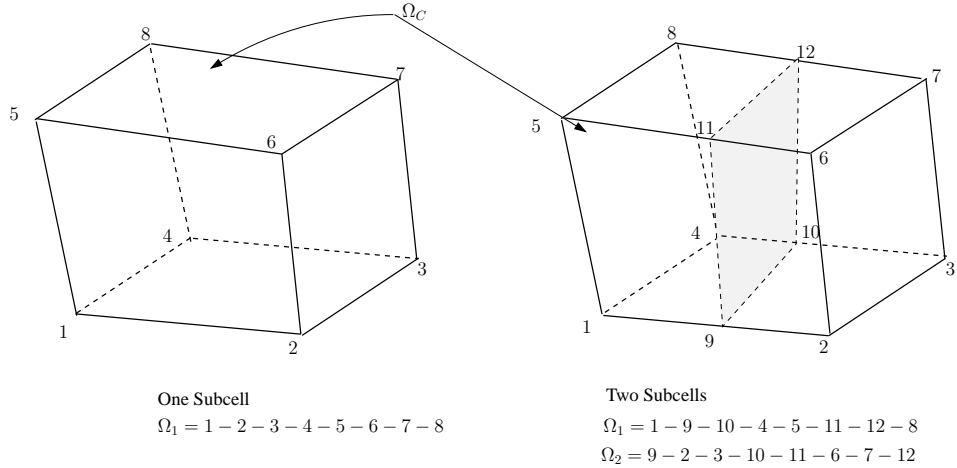

When the strain smoothing is used over three-dimensional domains, the volume integral is transferred to a surface integral. This surface integral is to be performed over the polygonal surfaces that build up the polyhedron. As in the 2D case, the smoothed element stiffness matrix is the sum over the subcells of the contribution from each subcell\footnote{Note that in 3D, the subcell is a volume} (see \fref{fig:subcellRepre3d}), which is constant:
\begin{equation}
\tilde{\kk}^e = \sum\limits_{C=1}^{nc}\int_{\Omega_C} {\tilde
{\bigb}_C^{\textup T} \dd\tilde {\bigb}_C } d\Omega =\sum\limits_{C=1}^{nc}
{\tilde {\bigb}_C^T \dd\tilde {\bigb}_C } \int_{\Omega_C}d\Omega
=\sum\limits_{C=1}^{nc}{\tilde {\bigb}^T_C \dd\tilde {\bigb}_C V_C }
\label{eqn:stiffnessvar3d}
\end{equation}
The strain-displacement matrix $\tilde {\bigb}_C$ is constant over each $\Omega_C$ and is of the following form:
\begin{equation}
\tilde \bigb_{C}=\left[ {\begin{array}{*{20}c} \tilde {\bigb}_{C1} &
\tilde {\bigb}_{C2} & \tilde {\bigb}_{C3} & \cdots \tilde {\bigb}_{Cn}
\end{array}} \right] \label{eq:Btilde3dbis}
\end{equation}
where for all shape functions $I \in \{1,\dots,n\}$, the $6 \times 3$ submatrix $\tilde {\bigb}_{CI}$ represents the contribution to the strain displacement matrix associated with shape function $I$ and cell $C$ and writes :
\begin{align}
\renewcommand\arraystretch{2}
\forall I\in\{1,2,\dots,n\},\forall C \in \{1,2,\dots nc \}\tilde
{\bigb}_{CI} &=\frac{1}{V_C} \int_{S_C}\cn^T(\xx) N_I(\xx) \rmd S \nonumber \\ 
&=
\frac{1}{V_C} \int_{S_C} \left[
\begin{array}{*{20}c}
 { n_x(\xx)} & 0 & 0 \\
 0 & { n_y(\xx)} & 0\\
 0 & 0 & n_z(\xx) \\
  n_y(\xx) &  n_x(\xx) & 0 \\
 0 & n_z(\xx) & n_y(\xx) \\
 n_z(\xx) & 0 & n_x(\xx) 
\end{array}  \right] N_I(\xx) \rmd S
\label{eqn:Btilde3d}
\end{align}
As in 2D, due to the process of strain smoothing, only the shape function is involved in the calculation of the field gradient and hence the stiffness matrix. In this study, over an arbitrary polygonal surface, we employ Wacshspress interpolants~\cite{wachspress1971}. In computing the strain-displacement matrix given by \Eref{eqn:Btilde3d} and the stiffness matrix, only the shape functions associated with the polygonal surface contribute to the integral. To evaluate the integral in \Eref{eqn:Btilde3d}, two schemes are adopted~\cite{gaintalischi2014}: (a) nodal quadrature and (b) conforming interpolant quadrature.

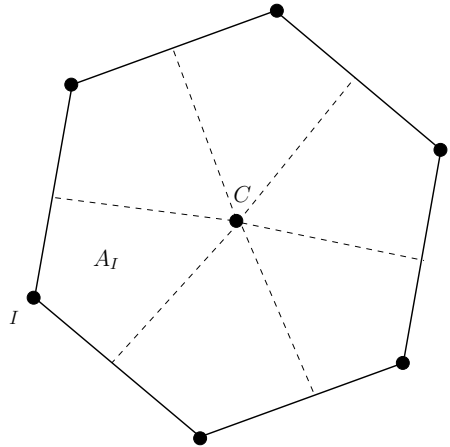
\begin{figure}[htpb]
\centering
\scalebox{0.6}{\input{./Figures/Nodalquad.pstex_t}}
\caption{Nodal quadrature, where $I$ is any node, $A_I$ is the area formed by the node, the centroid $C$ and the mid-points of the edges containing the node $I$.}
\label{fig:nodalquad}
\end{figure}

\paragraph{Nodal quadrature} In this case, the surface integral of the shape function $N_I$ over any face of the polyhedral element is given by:
\begin{equation}
\int\limits_{S_C} N_I (\xx) ~\rmd S = N_I(\xx_I) A_I = A_I
\end{equation}
where $A_I$ is the nodal weight of the node $I$, which is the area of the quadrilateral formed by the node, the centroid of the face and the mid-points of the edges containing the node. This is shown in \fref{fig:nodalquad}, however, this scheme is applicable only to the elements where star convexity is satisfied.

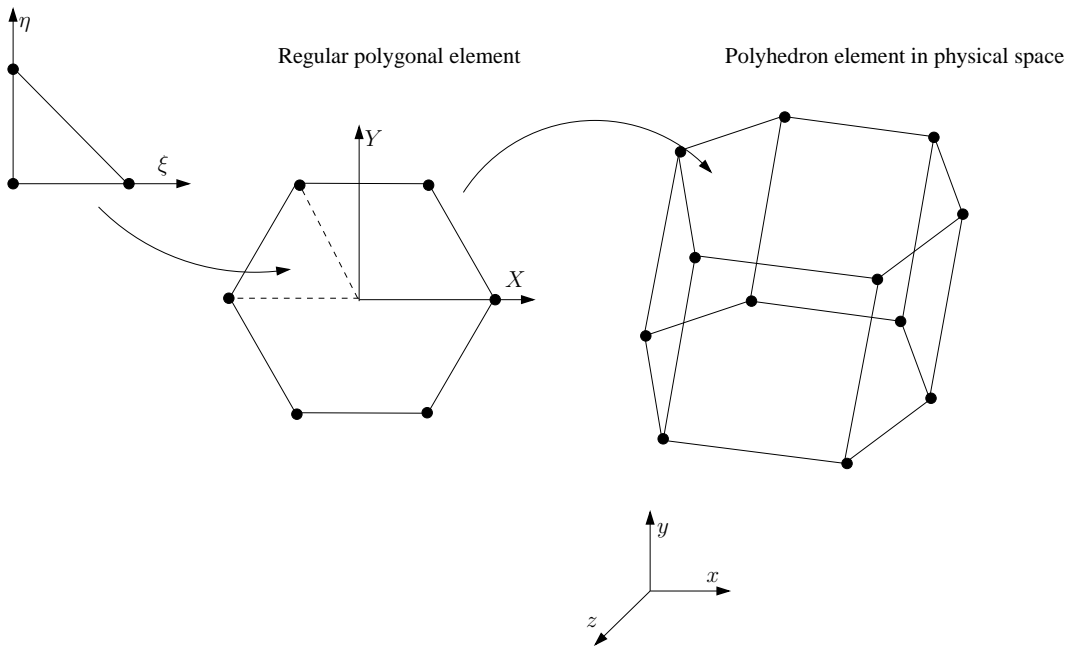
\begin{figure}[htpb]
\centering
\scalebox{0.6}{\input{./Figures/interpol.pstex_t}}
\caption{Surface integration scheme. The polygonal surface of a polyhedron is mapped onto a regular polygon. The regular polygon is then sub-divided into triangles and each of those sub-triangles is then mapped onto a standard triangle. Quadrature rules over the triangles are the used for numerical integration. 'Filled' circles denote the nodes.}
\label{fig:conformingsurfaceinte}
\end{figure}

\paragraph{Conforming Interpolant quadrature} In this case, to evaluate the surface integral of the shape function $N_I$, we adopt an interpolation scheme. For this purpose, the knowledge of the shape functions within the polygonal surface is required and in this study, we employ Wachspress interpolants~\cite{wachspress1971} over the polygonal surface. To integrate the terms in the strain-displacement matrix, the polygonal surface is mapped onto a regular polygon. The regular polygon is sub-divided into triangles and triangular quadrature rules over each triangle are employed to numerically integrate the terms in \Eref{eqn:Btilde3d} (see \fref{fig:conformingsurfaceinte}). This process involves a two level iso-parametric mapping of the surface and relies on the positivity of the Jacobian matrix involved in the transformation. The other possible approaches include: (a) complex mappings such as the Schwarz-Christoffel conformal mapping~\cite{balachandranrajagopal2008,natarajanbordas2009}; (b) adaptively weighted numerical integration scheme~\cite{thiagarajanshapiro2014}; (c) generalized Gaussian quadrature rules~\cite{mousavixiao2010} and (d) Guass-Green cubature~\cite{sommarivavianello2009}, to name a few.

%% file: Figures/sfemid.pstex_t
\begin{picture}(0,0)%
\includegraphics{./Figures/sfemid.pstex}%
\end{picture}%
\setlength{\unitlength}{4144sp}%
\begingroup\makeatletter\ifx\SetFigFont\undefined%
\gdef\SetFigFont#1#2#3#4#5{%
  \reset@font\fontsize{#1}{#2pt}%
  \fontfamily{#3}\fontseries{#4}\fontshape{#5}%
  \selectfont}%
\fi\endgroup%
\begin{picture}(11488,8020)(451,-7763)
\put(481,-4036){\makebox(0,0)[lb]{\smash{{\SetFigFont{14}{16.8}{\familydefault}{\mddefault}{\updefault}{\color[rgb]{0,0,0}(b)}%
}}}}
\put(2281, 14){\makebox(0,0)[lb]{\smash{{\SetFigFont{14}{16.8}{\familydefault}{\mddefault}{\updefault}{\color[rgb]{0,0,0}One Subcell}%
}}}}
\put(7681, -1){\makebox(0,0)[lb]{\smash{{\SetFigFont{14}{16.8}{\familydefault}{\mddefault}{\updefault}{\color[rgb]{0,0,0}n subcells}%
}}}}
\put(9466,-4621){\makebox(0,0)[lb]{\smash{{\SetFigFont{14}{16.8}{\familydefault}{\mddefault}{\updefault}{\color[rgb]{0,0,0}6 Triangular subcells}%
}}}}
\put(9586,-661){\makebox(0,0)[lb]{\smash{{\SetFigFont{14}{16.8}{\familydefault}{\mddefault}{\updefault}{\color[rgb]{0,0,0}4 quadrilateral subcells}%
}}}}
\put(466, 74){\makebox(0,0)[lb]{\smash{{\SetFigFont{14}{16.8}{\familydefault}{\mddefault}{\updefault}{\color[rgb]{0,0,0}(a)}%
}}}}
\end{picture}%

%% file: Figures/sfemid3d.pstex_t
\begin{picture}(0,0)%
\includegraphics{./Figures/sfemid3d.pstex}%
\end{picture}%
\setlength{\unitlength}{4144sp}%
\begingroup\makeatletter\ifx\SetFigFont\undefined%
\gdef\SetFigFont#1#2#3#4#5{%
  \reset@font\fontsize{#1}{#2pt}%
  \fontfamily{#3}\fontseries{#4}\fontshape{#5}%
  \selectfont}%
\fi\endgroup%
\begin{picture}(11138,5544)(1036,-5548)
\put(2866,-2866){\makebox(0,0)[lb]{\smash{{\SetFigFont{14}{16.8}{\familydefault}{\mddefault}{\updefault}{\color[rgb]{0,0,0}$4$}%
}}}}
\put(5701,-211){\makebox(0,0)[lb]{\smash{{\SetFigFont{14}{16.8}{\familydefault}{\mddefault}{\updefault}{\color[rgb]{0,0,0}$\Omega_C$}%
}}}}
\put(1636,-4006){\makebox(0,0)[lb]{\smash{{\SetFigFont{14}{16.8}{\familydefault}{\mddefault}{\updefault}{\color[rgb]{0,0,0}$1$}%
}}}}
\put(7231,-1576){\makebox(0,0)[lb]{\smash{{\SetFigFont{14}{16.8}{\familydefault}{\mddefault}{\updefault}{\color[rgb]{0,0,0}$5$}%
}}}}
\put(7831,-3946){\makebox(0,0)[lb]{\smash{{\SetFigFont{14}{16.8}{\familydefault}{\mddefault}{\updefault}{\color[rgb]{0,0,0}$1$}%
}}}}
\put(9286,-4111){\makebox(0,0)[lb]{\smash{{\SetFigFont{14}{16.8}{\familydefault}{\mddefault}{\updefault}{\color[rgb]{0,0,0}$9$}%
}}}}
\put(10636,-4186){\makebox(0,0)[lb]{\smash{{\SetFigFont{14}{16.8}{\familydefault}{\mddefault}{\updefault}{\color[rgb]{0,0,0}$2$}%
}}}}
\put(12136,-3256){\makebox(0,0)[lb]{\smash{{\SetFigFont{14}{16.8}{\familydefault}{\mddefault}{\updefault}{\color[rgb]{0,0,0}$3$}%
}}}}
\put(11971,-976){\makebox(0,0)[lb]{\smash{{\SetFigFont{14}{16.8}{\familydefault}{\mddefault}{\updefault}{\color[rgb]{0,0,0}$7$}%
}}}}
\put(10366,-691){\makebox(0,0)[lb]{\smash{{\SetFigFont{14}{16.8}{\familydefault}{\mddefault}{\updefault}{\color[rgb]{0,0,0}$12$}%
}}}}
\put(8836,-466){\makebox(0,0)[lb]{\smash{{\SetFigFont{14}{16.8}{\familydefault}{\mddefault}{\updefault}{\color[rgb]{0,0,0}$8$}%
}}}}
\put(9121,-2836){\makebox(0,0)[lb]{\smash{{\SetFigFont{14}{16.8}{\familydefault}{\mddefault}{\updefault}{\color[rgb]{0,0,0}$4$}%
}}}}
\put(9181,-1591){\makebox(0,0)[lb]{\smash{{\SetFigFont{14}{16.8}{\familydefault}{\mddefault}{\updefault}{\color[rgb]{0,0,0}$11$}%
}}}}
\put(4351,-4201){\makebox(0,0)[lb]{\smash{{\SetFigFont{14}{16.8}{\familydefault}{\mddefault}{\updefault}{\color[rgb]{0,0,0}$2$}%
}}}}
\put(5881,-3256){\makebox(0,0)[lb]{\smash{{\SetFigFont{14}{16.8}{\familydefault}{\mddefault}{\updefault}{\color[rgb]{0,0,0}$3$}%
}}}}
\put(1051,-1426){\makebox(0,0)[lb]{\smash{{\SetFigFont{14}{16.8}{\familydefault}{\mddefault}{\updefault}{\color[rgb]{0,0,0}$5$}%
}}}}
\put(2611,-481){\makebox(0,0)[lb]{\smash{{\SetFigFont{14}{16.8}{\familydefault}{\mddefault}{\updefault}{\color[rgb]{0,0,0}$8$}%
}}}}
\put(5626,-961){\makebox(0,0)[lb]{\smash{{\SetFigFont{14}{16.8}{\familydefault}{\mddefault}{\updefault}{\color[rgb]{0,0,0}$7$}%
}}}}
\put(4096,-1756){\makebox(0,0)[lb]{\smash{{\SetFigFont{14}{16.8}{\familydefault}{\mddefault}{\updefault}{\color[rgb]{0,0,0}$6$}%
}}}}
\put(8116,-4756){\makebox(0,0)[lb]{\smash{{\SetFigFont{14}{16.8}{\familydefault}{\mddefault}{\updefault}{\color[rgb]{0,0,0}Two Subcells}%
}}}}
\put(2806,-4801){\makebox(0,0)[lb]{\smash{{\SetFigFont{14}{16.8}{\familydefault}{\mddefault}{\updefault}{\color[rgb]{0,0,0}One Subcell}%
}}}}
\put(2806,-5116){\makebox(0,0)[lb]{\smash{{\SetFigFont{14}{16.8}{\familydefault}{\mddefault}{\updefault}{\color[rgb]{0,0,0}$\Omega_1 = 1-2-3-4-5-6-7-8$}%
}}}}
\put(8086,-5116){\makebox(0,0)[lb]{\smash{{\SetFigFont{14}{16.8}{\familydefault}{\mddefault}{\updefault}{\color[rgb]{0,0,0}$\Omega_1 = 1-9-10-4-5-11-12-8$}%
}}}}
\put(8056,-5461){\makebox(0,0)[lb]{\smash{{\SetFigFont{14}{16.8}{\familydefault}{\mddefault}{\updefault}{\color[rgb]{0,0,0}$\Omega_2 = 9-2-3-10-11-6-7-12$}%
}}}}
\put(10636,-1951){\makebox(0,0)[lb]{\smash{{\SetFigFont{14}{16.8}{\familydefault}{\mddefault}{\updefault}{\color[rgb]{0,0,0}$6$}%
}}}}
\put(10426,-2821){\makebox(0,0)[lb]{\smash{{\SetFigFont{14}{16.8}{\familydefault}{\mddefault}{\updefault}{\color[rgb]{0,0,0}$10$}%
}}}}
\end{picture}%

%% file: Figures/Nodalquad.pstex_t
\begin{picture}(0,0)%
\includegraphics{./Figures/Nodalquad.pstex}%
\end{picture}%
\setlength{\unitlength}{3947sp}%
\begingroup\makeatletter\ifx\SetFigFont\undefined%
\gdef\SetFigFont#1#2#3#4#5{%
  \reset@font\fontsize{#1}{#2pt}%
  \fontfamily{#3}\fontseries{#4}\fontshape{#5}%
  \selectfont}%
\fi\endgroup%
\begin{picture}(4560,4575)(3631,-5086)
\put(3646,-3826){\makebox(0,0)[lb]{\smash{{\SetFigFont{12}{14.4}{\familydefault}{\mddefault}{\updefault}{\color[rgb]{0,0,0}$I$}%
}}}}
\put(4516,-3226){\makebox(0,0)[lb]{\smash{{\SetFigFont{14}{16.8}{\familydefault}{\mddefault}{\updefault}{\color[rgb]{0,0,0}$A_I$}%
}}}}
\put(5971,-2566){\makebox(0,0)[lb]{\smash{{\SetFigFont{14}{16.8}{\familydefault}{\mddefault}{\updefault}{\color[rgb]{0,0,0}$C$}%
}}}}
\end{picture}%

%% file: Figures/interpol.pstex_t
\begin{picture}(0,0)%
\includegraphics{./Figures/interpol.pstex}%
\end{picture}%
\setlength{\unitlength}{3947sp}%
\begingroup\makeatletter\ifx\SetFigFont\undefined%
\gdef\SetFigFont#1#2#3#4#5{%
  \reset@font\fontsize{#1}{#2pt}%
  \fontfamily{#3}\fontseries{#4}\fontshape{#5}%
  \selectfont}%
\fi\endgroup%
\begin{picture}(11226,7257)(1139,-6958)
\put(8581,-916){\makebox(0,0)[lb]{\smash{{\SetFigFont{14}{16.8}{\familydefault}{\mddefault}{\updefault}{\color[rgb]{0,0,0}Polyhedron element in physical space}%
}}}}
\put(8386,-6286){\makebox(0,0)[lb]{\smash{{\SetFigFont{14}{16.8}{\familydefault}{\mddefault}{\updefault}{\color[rgb]{0,0,0}$x$}%
}}}}
\put(7876,-5776){\makebox(0,0)[lb]{\smash{{\SetFigFont{14}{16.8}{\familydefault}{\mddefault}{\updefault}{\color[rgb]{0,0,0}$y$}%
}}}}
\put(6301,-3226){\makebox(0,0)[lb]{\smash{{\SetFigFont{14}{16.8}{\familydefault}{\mddefault}{\updefault}{\color[rgb]{0,0,0}$X$}%
}}}}
\put(4846,-1771){\makebox(0,0)[lb]{\smash{{\SetFigFont{14}{16.8}{\familydefault}{\mddefault}{\updefault}{\color[rgb]{0,0,0}$Y$}%
}}}}
\put(2701,-2026){\makebox(0,0)[lb]{\smash{{\SetFigFont{14}{16.8}{\familydefault}{\mddefault}{\updefault}{\color[rgb]{0,0,0}$\xi$}%
}}}}
\put(1261,-541){\makebox(0,0)[lb]{\smash{{\SetFigFont{14}{16.8}{\familydefault}{\mddefault}{\updefault}{\color[rgb]{0,0,0}$\eta$}%
}}}}
\put(7141,-6751){\makebox(0,0)[lb]{\smash{{\SetFigFont{14}{16.8}{\familydefault}{\mddefault}{\updefault}{\color[rgb]{0,0,0}$z$}%
}}}}
\put(1156,104){\makebox(0,0)[lb]{\smash{{\SetFigFont{14}{16.8}{\familydefault}{\mddefault}{\updefault}{\color[rgb]{0,0,0}Standard triangular element}%
}}}}
\put(3946,-916){\makebox(0,0)[lb]{\smash{{\SetFigFont{14}{16.8}{\familydefault}{\mddefault}{\updefault}{\color[rgb]{0,0,0}Regular polygonal element}%
}}}}
\end{picture}%

%% file: vemintro.tex
The virtual element method can be seen as a generalization of the finite element method to arbitrary polygons and polyhedra. The VEM does not require quadrature formulae to compute the stiffness matrix nor an expression for the basis functions. The explicit computation of the basis functions is actually not needed and this is the reason of the word `Virtual' in the VEM~\cite{veigabrezzi2013}. The important ingredient is the operator $\mathbf{\Pi}^\nabla$ that relates to the bilinear form of the problem. Once this is known, the local element stiffness matrix can be computed. It is beyond the scope of this paper to discuss the details of the VEM as applied to scalar and elasticity problems. Only important equations pertaining to the computation of the stiffness matrix associated to an arbitrary polygon/polyhedron is given in this section. Interested readers are referred to  the work of Beir\~ao Da Veiga \textit{et al.,}~\cite{veigabrezzi2013} for scalar problems and Beir\~ao Da Veiga \textit{et al.,}~\cite{veigabrezzi2013a} and Gain \textit{et al.,}~\cite{gaintalischi2014} for three dimensional elasticity, where the method is dealt with in great detail. In this section, we only present the final expression to compute the stiffness matrix by the VEM in to three dimensional linear elasticity. For more detailed derivation and discussion, interested readers are referred to the work of Gain \textit{et al.,}~\cite{gaintalischi2014}. The expression for the stiffness matrix can be written as:
\begin{equation}
\kk^E_h  = \underbrace{|E| \mathbf{W}_C \dd \mathbf{W}_C^{\rm T}}_{\rm Consistency~term = \kk_{\rm const}} + \underbrace{(\mathbf{I} - \mathbf{P}_p)^{\rm T} \mathbf{S}^E (\mathbf{I} - \mathbf{P}_p)}_{\rm Stability~term = \kk_{\rm stab}}
\label{eqn:vem3d}
\end{equation}
where $|E|$ is the measure of the polygon and
\begin{equation}
\mathbf{P}_p = \mathbf{P}_R + \mathbf{P}_C
\end{equation}
and
\begin{align}
\mathbf{P}_R &= \mathbf{N}_R \mathbf{W}_R^{\rm T} \nonumber \\
\mathbf{P}_C &= \mathbf{N}_C \mathbf{W}_C^{\rm T} 
\end{align}
The block $3I-2 \colon 3I$ rows of $\mathbf{N}_R$ and $\mathbf{N}_C $ are expressed as:
\begin{align}
\mathbf{N}_R(3I-2 \colon 3) &=\left[ \begin{array}{rrrrrr} 1 & 0 & 0 & (\xx_I-\overline{\xx})_2 & 0 & (\xx_I-\overline{\xx})_3 \\ 0 & 1 & 0 & -(\xx_I-\overline{\xx})_1 & (\xx_I-\overline{\xx})_3 & 0 \\ 0 & 0 & 1 & 0 & -(\xx_I-\overline{\xx})_2 & (\xx_I-\overline{\xx})_1 \end{array} \right] \nonumber \\
\mathbf{N}_C (3I-2 \colon 3) &= \left[ \begin{array}{rrrrrr} (\xx_I-\overline{\xx})_1 & 0 & 0 & (\xx_I-\overline{\xx})_2 & 0 & (\xx_I-\overline{\xx})_3 \\ 0 & (\xx_I-\overline{\xx})_2 & (\xx_I-\overline{\xx})_1 & (\xx_I-\overline{\xx})_3 & 0 \\ 0 & 0 & (\xx_I-\overline{\xx})_3 & 0 & (\xx_I-\overline{\xx})_2 & (\xx_I-\overline{\xx})_1 \end{array} \right]
\end{align}
where $\xx_I$ is the coordinate of the node and $\overline{\xx}$ is the polyhedron centroid. The block $3I-2 \colon 3I$ rows of $\mathbf{W}_R$ and $\mathbf{W}_C$ are expressed as:
\begin{align}
\mathbf{W}_R(3I-2 \colon 3) &= \left[ \begin{array}{rrrrrr} 1/n & 0 & 0 & (\mathbf{q}_I)_2 & 0 & -(\mathbf{q}_I)_3 \\ 0 & 1/n & 0 & -(\mathbf{q}_I)_1 & (\mathbf{q}_I)_3 & 0 \\ 0 & 0 & 1/n & 0 -(\mathbf{q}_I)_2 & (\mathbf{q}_I)_1 \end{array} \right] \nonumber \\
\mathbf{W}_C(3I-2 \colon 3) &= \left[ \begin{array}{rrrrrr} 2(\mathbf{q}_I)_1 & 0 & 0 & (\mathbf{q}_I)_2 & 0 & (\mathbf{q}_I)_3 \\ 0 & 2(\mathbf{q}_I)_2 & 0 & (\mathbf{q}_I)_1 & (\mathbf{q}_I)_3 & 0 \\ 0 & 0 & 2(\mathbf{q}_I)_3 & 0 & (\mathbf{q}_I)_2 & (\mathbf{q}_I)_1 \end{array} \right]
\end{align}
where the subscript indicates the component of the associated vector and 
\begin{equation*}
\mathbf{q}_I = \frac{1}{2|E|} \int\limits_{\partial E} N_I \mathbf{n}~\rmd \Gamma
\end{equation*}
where $\mathbf{S}^E = \alpha \mathbf{I}$ and $\alpha = \alpha^\ast \textup{trace}( |E| \mathbf{W}_C \dd \mathbf{W}_C^{\rm T})$ is a scaling coefficient and $\mathbf{n}$ is the unit outward normal vector. It can be seen that the computation of the stiffness matrix involves computing the matrices $\mathbf{N}_R, \mathbf{N}_C, \mathbf{W}_R$ and $\mathbf{W}_C$. The calculation of the matrices $\mathbf{W}_R$ and $\mathbf{W}_C$ involves computing the surface integral of the basis functions. This can be computed by employing one of the techniques discussed in the previous section.

%% file: twoDexamples.tex
\subsubsection{Cantilever beam}
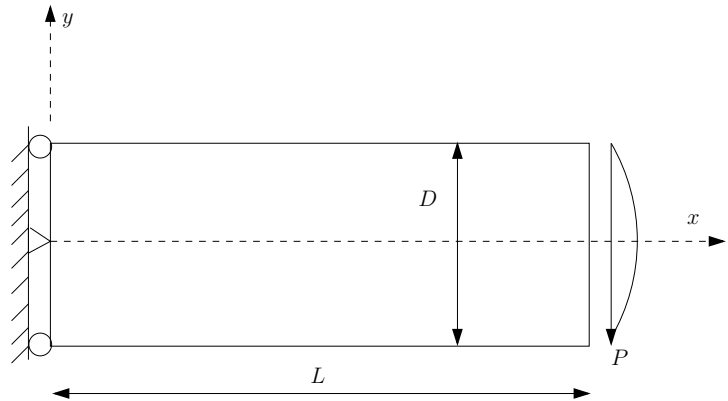
\begin{figure}[htpb]
\centering
\scalebox{0.7}{\input{./Figures/canti.pstex_t}}
\caption{Cantilever beam: Geometry and boundary conditions.}
\label{fig:cantmesh}
\end{figure}

\begin{figure}[htpb]
\centering
\includegraphics[scale=0.7]{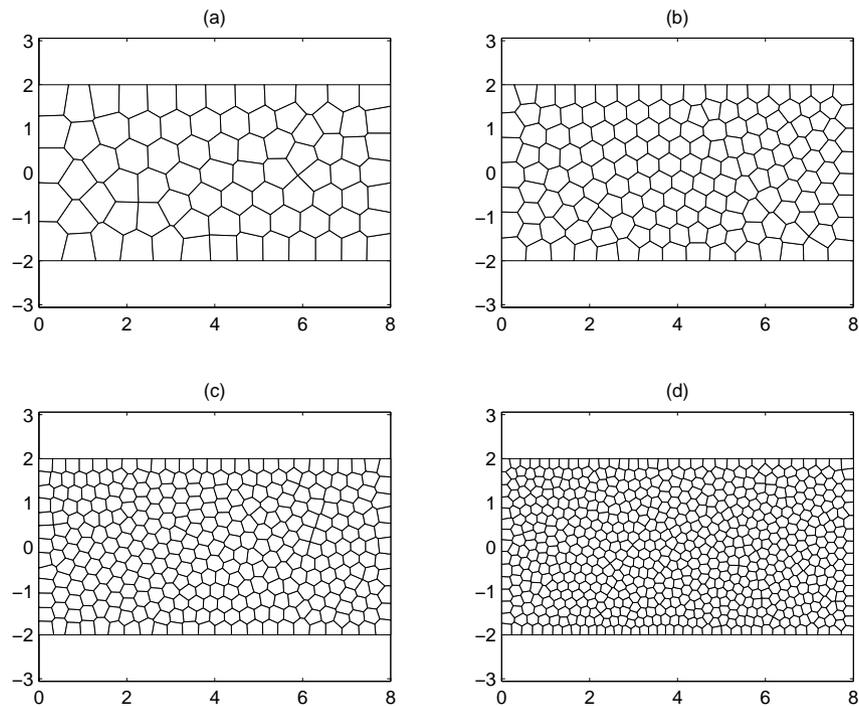}
\caption{Cantilever beam: Typical polygonal mesh employed in this study.}
\label{fig:cantileverfig}
\end{figure}

\begin{figure}[htpb]
\centering
\newlength\figureheight 
\newlength\figurewidth 
\setlength\figureheight{8cm} 
\setlength\figurewidth{10cm}
\input{./Figures/cantibeamInfluAlpha.tikz}
\caption{Bending of thick cantilever beam: Influence of $\alpha^\ast$ for mesh (d) (c.f. \fref{fig:cantileverfig}) }
\label{fig:cantiConveResults_alpha}
\end{figure}
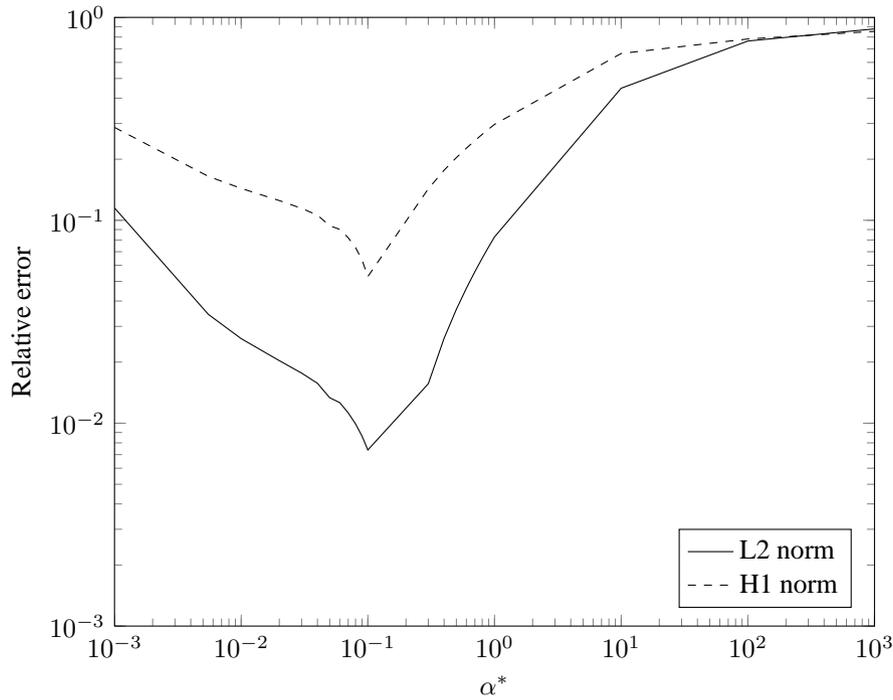

A two-dimensional cantilever beam subjected to a parabolic shear load at the free end is examined as shown in \fref{fig:cantmesh}. The geometry is: length $L=$ 8m, height $D=$ 4m. The material properties are: Young's modulus, $E=$ 3e$^7$ N/m$^2$, Poisson's ratio $\nu=$ 0.3 and the parabolic shear force $P=$ 250 N. The exact solution for displacements are given by:
\begin{align}
u(x,y) &= \frac{P y}{6 \overline{E}I} \left[ (9L-3x)x + (2+\overline{\nu}) \left( y^2 - \frac{D^2}{4} \right) \right] \nonumber \\
v(x,y) &= -\frac{P}{6 \overline{E}I} \left[ 3\overline{\nu}y^2(L-x) + (4+5\overline{\nu}) \frac{D^2x}{4} + (3L-x)x^2 \right]
\label{eqn:cantisolution}
\end{align}
where $I = D^3/12$ is the moment of inertia, $\overline{E} = E$, $\overline{\nu} = \nu$ and $\overline{E} = E/(1-\nu^2)$, $\overline{\nu} = \nu/(1-\nu)$ for plane stress and plane strain, respectively. The domain is discretized with two different mesh types: (a) structured quadrilateral elements (8$\times$4, 16$\times$8, 32$\times$16, 64$\times$32) and (b) polygonal elements. \fref{fig:cantileverfig} shows a sample polygonal mesh used for this study. Before demonstrating the convergence and the accuracy of the SFEM with stabilization, we investigate the influence of the scaling parameter $\alpha^\ast$. \fref{fig:cantiConveResults_alpha} shows the influence of the scaling parameter $\alpha^\ast$ on the relative error in the $L^2$ and $H^1$ norm. It is observed that the relative error attains a minimum value for $\alpha^\ast=$ 0.1. 

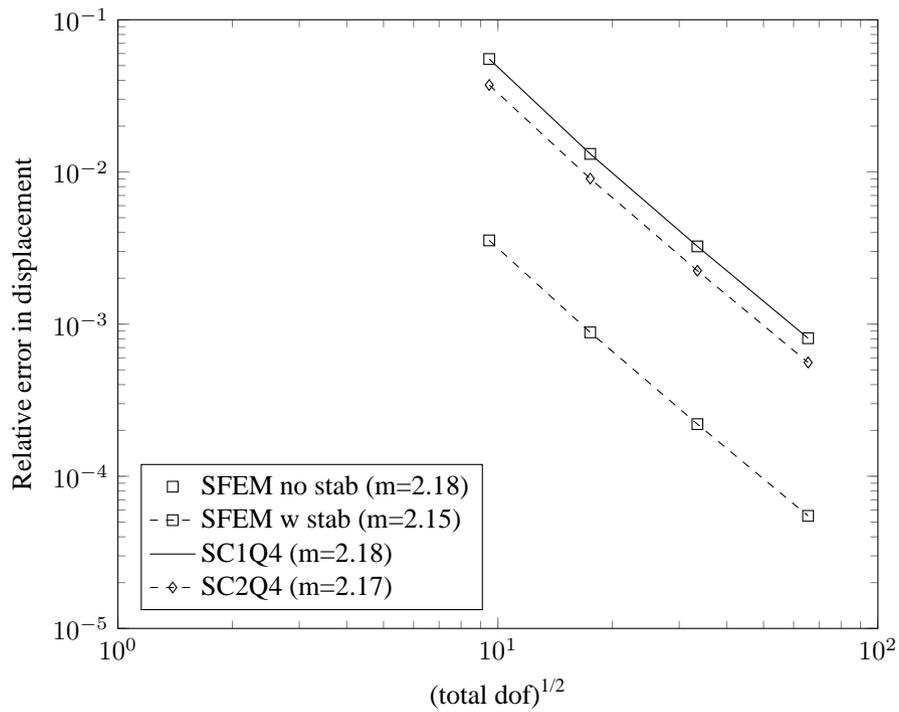
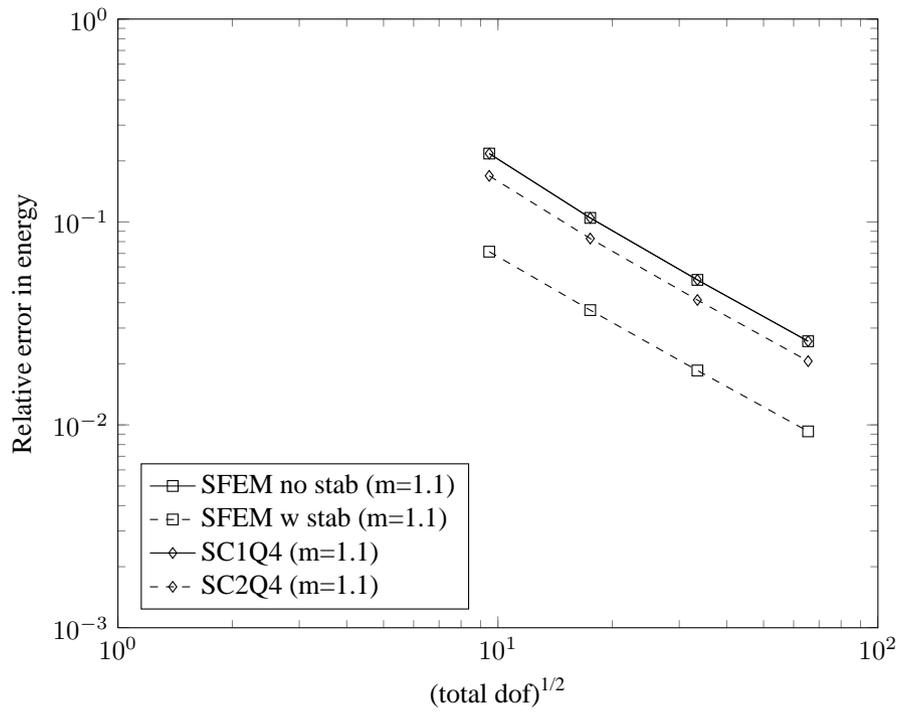
\begin{figure}
\centering
\setlength\figureheight{8cm} 
\setlength\figurewidth{10cm}
\subfigure[]{\input{./Figures/cantibeamL2.tikz}}
\subfigure[]{\input{./Figures/cantibeamH1.tikz}}
\caption{Bending of thick cantilever beam: Convergence results for (a) the relative error in the displacement norm $(L^2)$ and (b) the relative error in the energy norm. The rate of convergence is also shown, where $m$ is the average slope. The domain is discretized with Q4 elements.}
\label{fig:cantiConveResults}
\end{figure}

The numerical convergence of the relative error in the displacement norm and the relative error in the energy norm is shown in \fref{fig:cantiConveResults} for structured quadrilateral elements. The problem is solved with conventional SFEM with one (SC1Q4) and two subcells (SC2Q4) and with the proposed SFEM with and without stabilization. It is observed that the SC1Q4 and SFEM with no stabilization yield similar results, as seen in the scalar example in Section \ref{sfemvsvem}. However, the SC2Q4 and SFEM with stabilization yield different results. This can be attributed to the choice of the scaling parameter $\alpha^\ast$. It can be seen that with mesh refinement, all the approaches converge with optimal rate. For a proper choice of $\alpha^\ast$, it is seen that SFEM with stabilization yields more accurate results than the SC2Q4. 

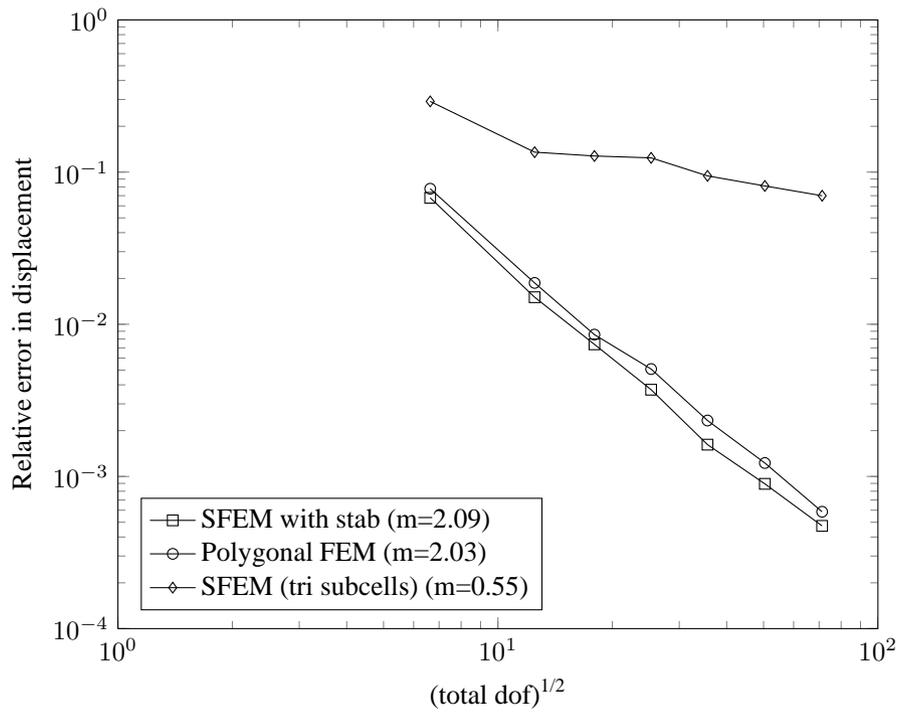
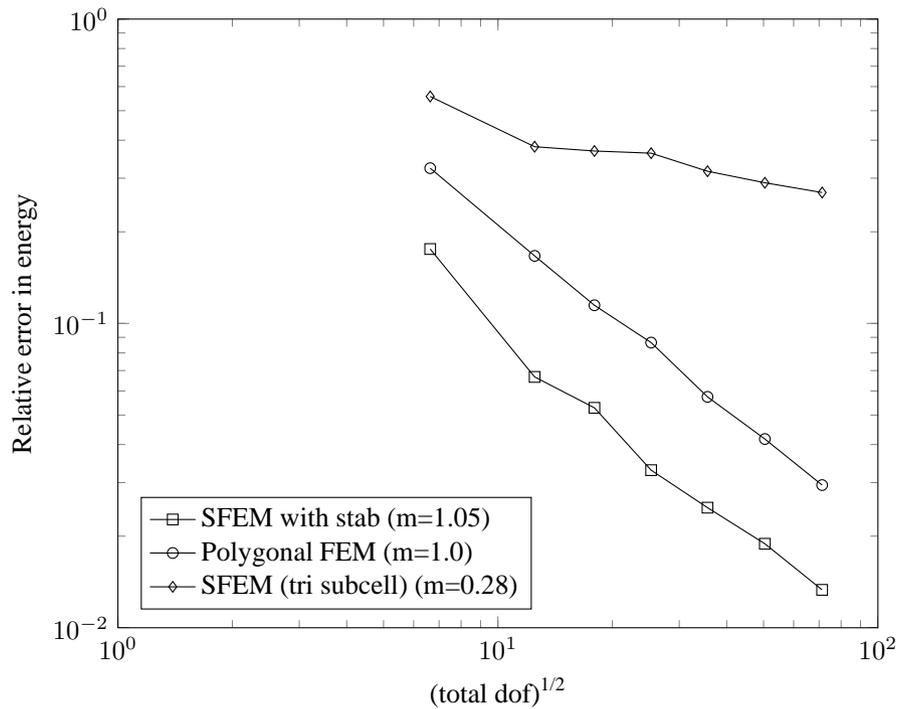
\begin{figure}
\centering
\setlength\figureheight{8cm} 
\setlength\figurewidth{10cm}
\subfigure[]{\input{./Figures/cantibeamL2Poly.tikz}}
\subfigure[]{\input{./Figures/cantibeamH1Poly.tikz}}
\caption{Bending of thick cantilever beam: Convergence results for (a) the relative error in the displacement norm $(L^2)$ and (b) the relative error in the energy norm. The rate of convergence is also shown, where $m$ is the average slope. The domain is discretized with arbitrary polygonal elements.}
\label{fig:cantiConveResultsPoly}
\end{figure}

\fref{fig:cantiConveResultsPoly} shows the convergence of the displacement and the energy norm with mesh refinement. In this case, one subcell per polygonal element cannot be used, as this would lead to spurious energy modes~\cite{dailiu2007}. Hence, in this study, for the conventional SFEM and the polygonal FEM, we sub-triangulate the polygonal element and integrate over each sub-triangle. In the case of the proposed approach, we employ one subcell and add the stabilization term (\Eref{eqn:stabterm}). It can be observed that the proposed SFEM with stabilization yields more accurate results compared to the conventional SFEM with triangulation. It is also seen that the polygonal FEM and the proposed approach converges at an optimal rate for both the displacement norm and the energy norm.

\subsubsection{Infinite plate with a circular hole}
\begin{figure}[htpb]
\centering
\scalebox{0.9}{\input{./Figures/platehole.pstex_t}}
\caption{Infinite plate with a circular hole.}
\label{fig:twocrkhole}
\end{figure}
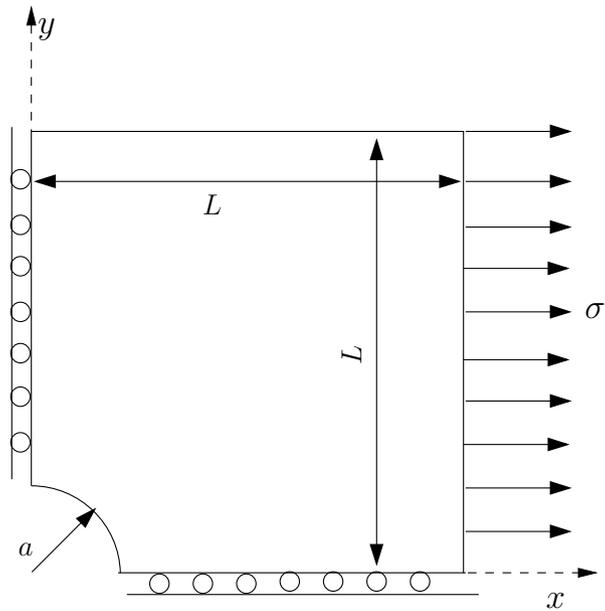

\begin{figure}[htpb]
\centering
\includegraphics[scale=0.75]{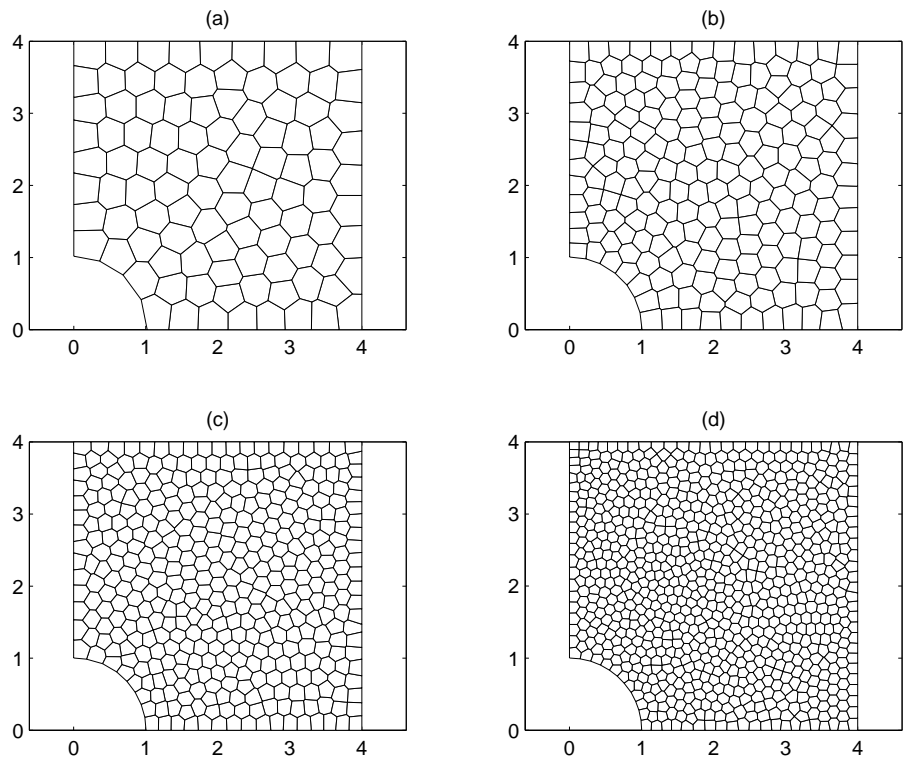}
\caption{Plate with a circular hole: domain discretized with polygonal elements: (a) 100 elements; (b) 200 elements; (c) 400 elements and (d) 800 elements.}
\label{fig:phwmesh}
\end{figure}

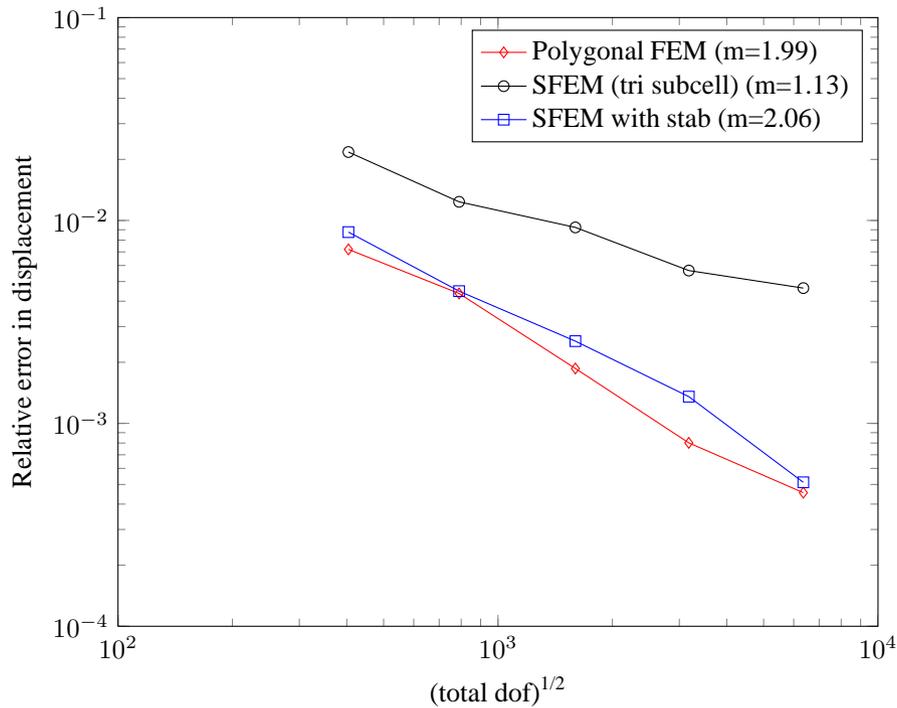
\begin{figure}
\centering
\setlength\figureheight{8cm} 
\setlength\figurewidth{10cm}
\input{./Figures/pwHoledError.tikz}
\caption{Infinite plate with a circular hole: Convergence results for the relative error in the displacement norm $(L^2)$. The rate of convergence is also shown, where $m$ is the average slope.}
\label{fig:plateHoleConveResults}
\end{figure}

In this example, consider an infinite plate with a traction free hole under uniaxial tension $(\sigma=$1$)$ along the $x-$axis~\fref{fig:twocrkhole}. The exact solution of the principal stresses in polar coordinates $(r,\theta)$ is given by:
\begin{eqnarray}
\sigma_{11}(r,\theta) &= 1 - \frac{a^2}{r^2} \left( \frac{3}{2} (\cos 2\theta + \cos 4\theta) \right) + \frac{3a^4}{2r^4} \cos 4\theta \nonumber \\
\sigma_{22}(r,\theta) &= -\frac{a^2}{r^2} \left( \frac{1}{2}(\cos 2\theta - \cos 4\theta) \right) - \frac{3a^4}{2r^4} \cos 4\theta \nonumber \\
\sigma_{12}(r,\theta) &= -\frac{a^2}{r^2} \left( \frac{1}{2}(\sin 2\theta + \sin 4\theta) \right) + \frac{3a^4}{2r^4} \sin 4\theta
\end{eqnarray}
where $a$ is the radius of the hole. Owing to symmetry, only one quarter of the plate is modeled. \fref{fig:phwmesh} shows a typical polygonal mesh used for the study. The material properties are: Young's modulus $E=$ 10$^5$ and Poisson's ratio $\nu=$ 0.3. In this example, analytical tractions are applied on the boundary. The domain is discretized with polygonal elements and along each edge of each polygon, the shape function is linear. The convergence rate in terms of the displacement norm is shown in \fref{fig:plateHoleConveResults}. It is observed that the proposed method yields more accurate results compared to the conventional SFEM with triangular subcells. For the present study, the scaling coefficient $\alpha^\ast$ is taken as 0.1.

\subsubsection{L-shaped domain under mode I loading}

\begin{figure}[htpb]
\centering
\includegraphics[scale=1.5]{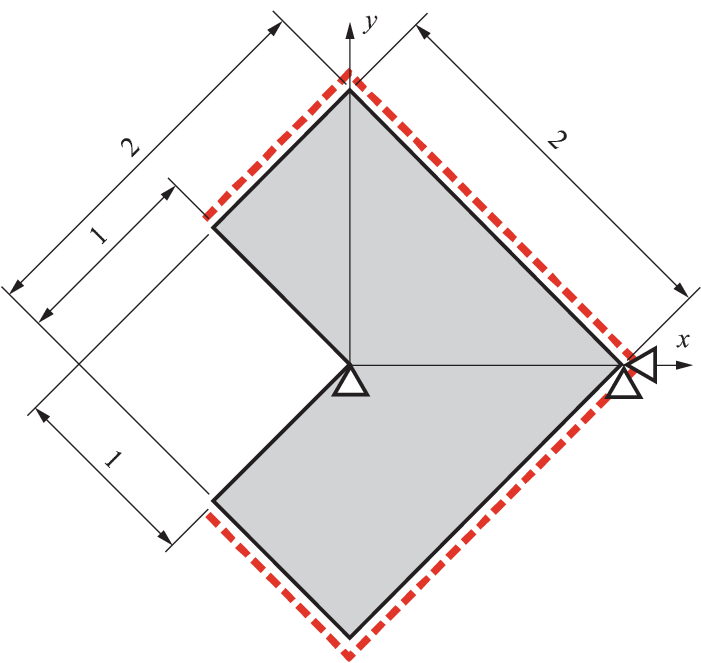}
\caption{L-shaped domain.}
\label{fig:LShape}
\end{figure}

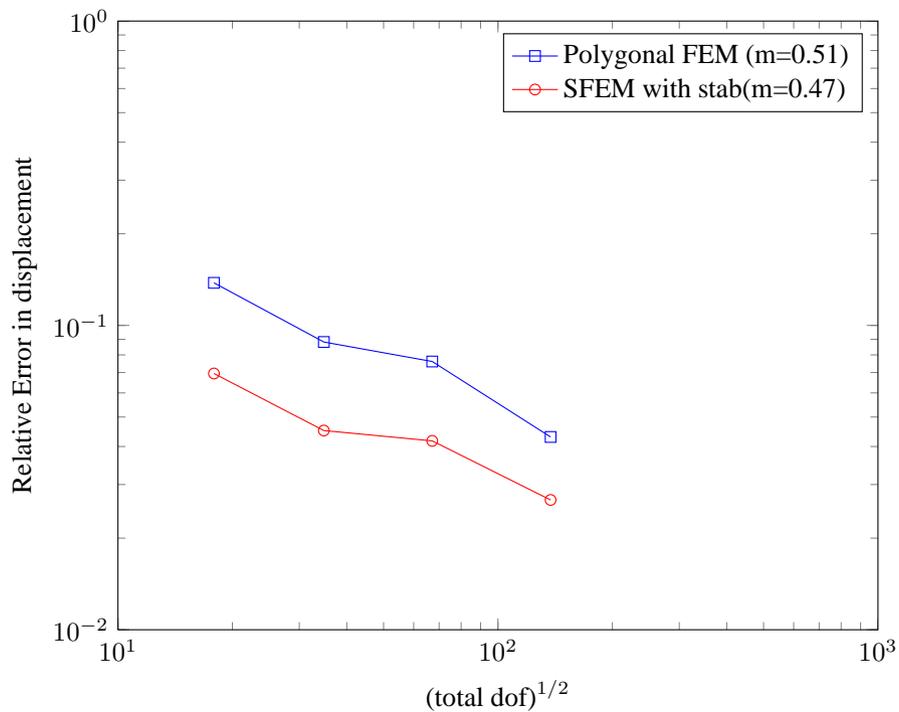
\begin{figure}
\centering
\setlength\figureheight{8cm} 
\setlength\figurewidth{10cm}
\input{./Figures/LshapeError.tikz}
\caption{L-shaped domain under mode I loading: Convergence results for the relative error in the displacement norm $(L^2)$. The rate of convergence is also shown, where $m$ is the average slope.}
\label{fig:Lshapeconveresults}
\end{figure}

In this example, consider the singular problem of a finite portion of an infinite domain with a reentrant corner. 
The model is loaded on part of the boundary, which is indicated by discontinuous thick lines in~Figure~\ref{fig:LShape}. The tractions correspond to the first terms of the asymptotic expansion that describes the exact solution under mixed mode loading conditions around the singular vertex. The exact displacement and stress fields for this singular elasticity problem can be found in \cite{szabobabuvska1991}. Exact values of the generalised stress intensity factors (GSIF) \cite{szabobabuvska1991} under mode I were taken as $K_{\rm I} = 1$  and $ K_{\rm II} = 0$. The material parameters are Young's modulus $E~=~1000$, and Poisson's ratio $\nu=$0.3 and the domain is discretized with polygonal elements. The problem is solved by the conventional polygonal FEM and the SFEM with stabilization. The convergence of the relative error in displacement with mesh refinement is shown in \fref{fig:Lshapeconveresults}. It is observed that both the approaches converge with mesh refinement. However, the SFEM with stabilization yields more accurate results. It is noted that since the domain has a reentrant corner, the optimal convergence rate is not achieved.

%% file: Figures/canti.pstex_t
\begin{picture}(0,0)%
\includegraphics{./Figures/canti.pstex}%
\end{picture}%
\setlength{\unitlength}{3947sp}%
\begingroup\makeatletter\ifx\SetFigFont\undefined%
\gdef\SetFigFont#1#2#3#4#5{%
  \reset@font\fontsize{#1}{#2pt}%
  \fontfamily{#3}\fontseries{#4}\fontshape{#5}%
  \selectfont}%
\fi\endgroup%
\begin{picture}(6369,3517)(2629,-3830)
\put(3091,-496){\makebox(0,0)[lb]{\smash{{\SetFigFont{12}{14.4}{\familydefault}{\mddefault}{\updefault}{\color[rgb]{0,0,0}$y$}%
}}}}
\put(8641,-2281){\makebox(0,0)[lb]{\smash{{\SetFigFont{12}{14.4}{\familydefault}{\mddefault}{\updefault}{\color[rgb]{0,0,0}$x$}%
}}}}
\put(6256,-2116){\makebox(0,0)[lb]{\smash{{\SetFigFont{12}{14.4}{\familydefault}{\mddefault}{\updefault}{\color[rgb]{0,0,0}$D$}%
}}}}
\put(5296,-3691){\makebox(0,0)[lb]{\smash{{\SetFigFont{12}{14.4}{\familydefault}{\mddefault}{\updefault}{\color[rgb]{0,0,0}$L$}%
}}}}
\put(7966,-3526){\makebox(0,0)[lb]{\smash{{\SetFigFont{12}{14.4}{\familydefault}{\mddefault}{\updefault}{\color[rgb]{0,0,0}$P$}%
}}}}
\end{picture}%

%% file: Figures/cantibeamInfluAlpha.tikz
%
%
\begin{tikzpicture}

\begin{axis}[%
width=\figurewidth,
height=\figureheight,
scale only axis,
xmode=log,
xmin=0.001,
xmax=1000,
xminorticks=true,
xlabel={$\alpha{}^\ast$},
ymode=log,
ymin=0.001,
ymax=1,
yminorticks=true,
ylabel={Relative error},
legend style={draw=black,fill=white,legend cell align=left},
legend pos = south east
]
\addplot [color=black,solid]
  table[row sep=crcr]{0.001	0.115111277138281\\
0.0055	0.034401911689718\\
0.01	0.026108377973161\\
0.02	0.020337816418117\\
0.03	0.01765756093918\\
0.04	0.015730081960731\\
0.05	0.01335151674551\\
0.06	0.012628160624641\\
0.07	0.011243065751459\\
0.08	0.009916791843523\\
0.09	0.008632605233314\\
0.1	0.007381089751062\\
0.3	0.015643426063211\\
0.4	0.026221820145679\\
0.5	0.036419059558411\\
0.6	0.046275084440062\\
0.7	0.055819724799878\\
0.8	0.065076303717194\\
0.9	0.07406395278225\\
1	0.082798914887165\\
10	0.447636252069992\\
100	0.764800083444511\\
1000	0.877692646672861\\
};
\addlegendentry{L2 norm};

\addplot [color=black,dashed]
  table[row sep=crcr]{0.001	0.286970609162176\\
0.0055	0.164667995400524\\
0.01	0.143638892142534\\
0.02	0.125056107934909\\
0.03	0.114521108495597\\
0.04	0.105976649622107\\
0.05	0.09411162262053\\
0.06	0.090156769866613\\
0.07	0.082017688120625\\
0.08	0.073343341890936\\
0.09	0.063793314832165\\
0.1	0.052805323509268\\
0.3	0.142969566854981\\
0.4	0.176534190460942\\
0.5	0.203656299210617\\
0.6	0.226769679100027\\
0.7	0.24707074756891\\
0.8	0.265254678482294\\
0.9	0.281768725724435\\
1	0.296921259029388\\
10	0.664027685555605\\
100	0.782812255817949\\
1000	0.85336855294175\\
};
\addlegendentry{H1 norm};

\end{axis}
\end{tikzpicture}%

%% file: Figures/cantibeamL2.tikz
%
%
\begin{tikzpicture}

\begin{axis}[%
width=\figurewidth,
height=\figureheight,
scale only axis,
xmode=log,
xmin=1,
xmax=100,
xminorticks=true,
xlabel={$\text{(total dof)}^{\text{1/2}}$},
ymode=log,
ymin=1e-05,
ymax=0.1,
yminorticks=true,
ylabel={Relative error in displacement},
legend style={draw=black,fill=white,legend cell align=left},
legend pos = south west
]
\addplot [color=black,only marks,mark=square,mark options={solid}]
  table[row sep=crcr]{9.48683298050514	0.055216919839692\\
17.4928556845359	0.013137462828774\\
33.4962684488885	0.003242418350484\\
65.4980915752513	0.000807515857044248\\
};
\addlegendentry{SFEM no stab (m=2.18)};

\addplot [color=black,dashed,mark=square,mark options={solid}]
  table[row sep=crcr]{9.48683298050514	0.003549087887556\\
17.4928556845359	0.000882212093802937\\
33.4962684488885	0.00022014088460527\\
65.4980915752513	5.49359917698314e-05\\
};
\addlegendentry{SFEM w stab (m=2.15)};

\addplot [color=black,solid]
  table[row sep=crcr]{9.48683298050514	0.055216919839692\\
17.4928556845359	0.013137462828774\\
33.4962684488885	0.003242418350484\\
65.4980915752513	0.000807515857044248\\
};
\addlegendentry{SC1Q4 (m=2.18)};

\addplot [color=black,dashed,mark=diamond,mark options={solid}]
  table[row sep=crcr]{9.48683298050514	0.037250558105718\\
17.4928556845359	0.009039082899026\\
33.4962684488885	0.002242765339262\\
65.4980915752513	0.000559488547308263\\
};
\addlegendentry{SC2Q4 (m=2.17)};

\end{axis}
\end{tikzpicture}%

%% file: Figures/cantibeamH1.tikz
%
%
\begin{tikzpicture}

\begin{axis}[%
width=\figurewidth,
height=\figureheight,
scale only axis,
xmode=log,
xmin=1,
xmax=100,
xminorticks=true,
xlabel={$\text{(total dof)}^{\text{1/2}}$},
ymode=log,
ymin=0.001,
ymax=1,
yminorticks=true,
ylabel={Relative error in energy},
legend style={draw=black,fill=white,legend cell align=left},
legend pos = south west
]
\addplot [color=black,solid,mark=square,mark options={solid}]
  table[row sep=crcr]{9.48683298050514	0.217324436863771\\
17.4928556845359	0.10473118027895\\
33.4962684488885	0.051826759141135\\
65.4980915752513	0.0258366048591\\
};
\addlegendentry{SFEM no stab (m=1.1)};

\addplot [color=black,dashed,mark=square,mark options={solid}]
  table[row sep=crcr]{9.48683298050514	0.071396284619404\\
17.4928556845359	0.036784016352368\\
33.4962684488885	0.01853787388742\\
65.4980915752513	0.009288354299887\\
};
\addlegendentry{SFEM w stab (m=1.1)};

\addplot [color=black,solid,mark=diamond,mark options={solid}]
  table[row sep=crcr]{9.48683298050514	0.217324436863771\\
17.4928556845359	0.10473118027895\\
33.4962684488885	0.051826759141135\\
65.4980915752513	0.0258366048591\\
};
\addlegendentry{SC1Q4 (m=1.1)};

\addplot [color=black,dashed,mark=diamond,mark options={solid}]
  table[row sep=crcr]{9.48683298050514	0.168816562386009\\
17.4928556845359	0.082858336074459\\
33.4962684488885	0.041255137476257\\
65.4980915752513	0.02060886319596\\
};
\addlegendentry{SC2Q4 (m=1.1)};

\end{axis}
\end{tikzpicture}%

%% file: Figures/cantibeamL2Poly.tikz
%
%
\begin{tikzpicture}

\begin{axis}[%
width=\figurewidth,
height=\figureheight,
scale only axis,
xmode=log,
xmin=1,
xmax=100,
xminorticks=true,
xlabel={$\text{(total dof)}^{\text{1/2}}$},
ymode=log,
ymin=0.0001,
ymax=1,
yminorticks=true,
ylabel={Relative error in displacement},
legend style={draw=black,fill=white,legend cell align=left},
legend pos = south west
]
\addplot [color=black,solid,mark=square,mark options={solid}]
  table[row sep=crcr]{6.6332495807108	0.067783145276453\\
12.4899959967968	0.015059721116259\\
17.9443584449264	0.007381089751062\\
25.298221281347	0.003716265013615\\
35.6089876295297	0.001618281062587\\
50.3984126734166	0.000894159464689231\\
71.3021738799035	0.000473362302081062\\
};
\addlegendentry{SFEM with stab (m=2.09)};

\addplot [color=black,solid,mark=o,mark options={solid}]
  table[row sep=crcr]{6.6332495807108	0.077747621021312\\
12.4899959967968	0.018698601784431\\
17.9443584449264	0.00855163304121\\
25.298221281347	0.005079705086802\\
35.6089876295297	0.002331086339671\\
50.3984126734166	0.001225639044603\\
71.3021738799035	0.000586100836017659\\
};
\addlegendentry{Polygonal FEM (m=2.03)};

\addplot [color=black,solid,mark=diamond,mark options={solid}]
  table[row sep=crcr]{6.6332495807108	0.291431991964912\\
12.4899959967968	0.135201853056681\\
17.9443584449264	0.127599614196438\\
25.298221281347	0.12403678616964\\
35.6089876295297	0.094426827597117\\
50.3984126734166	0.081090438685846\\
71.3021738799035	0.069931146037355\\
};
\addlegendentry{SFEM (tri subcells) (m=0.55)};

\end{axis}
\end{tikzpicture}%

%% file: Figures/cantibeamH1Poly.tikz
%
%
\begin{tikzpicture}

\begin{axis}[%
width=\figurewidth,
height=\figureheight,
scale only axis,
xmode=log,
xmin=1,
xmax=100,
xminorticks=true,
xlabel={$\text{(total dof)}^{\text{1/2}}$},
ymode=log,
ymin=0.01,
ymax=1,
yminorticks=true,
ylabel={Relative error in energy},
legend style={draw=black,fill=white,legend cell align=left},
legend pos = south west
]
\addplot [color=black,solid,mark=square,mark options={solid}]
  table[row sep=crcr]{6.6332495807108	0.175573990747497\\
12.4899959967968	0.06671928170249\\
17.9443584449264	0.052805323509268\\
25.298221281347	0.032956024159045\\
35.6089876295297	0.024808402745654\\
50.3984126734166	0.018884459086623\\
71.3021738799035	0.013322352458685\\
};
\addlegendentry{SFEM with stab (m=1.05)};

\addplot [color=black,solid,mark=o,mark options={solid}]
  table[row sep=crcr]{6.6332495807108	0.323672955901606\\
12.4899959967968	0.166698087401175\\
17.9443584449264	0.114726583816981\\
25.298221281347	0.086497232140651\\
35.6089876295297	0.057365908665105\\
50.3984126734166	0.041706716349184\\
71.3021738799035	0.02943431168056\\
};
\addlegendentry{Polygonal FEM (m=1.0)};

\addplot [color=black,solid,mark=diamond,mark options={solid}]
  table[row sep=crcr]{6.6332495807108	0.556787643981744\\
12.4899959967968	0.380724526396742\\
17.9443584449264	0.368563098032213\\
25.298221281347	0.362591816450448\\
35.6089876295297	0.316265456557697\\
50.3984126734166	0.289845734737707\\
71.3021738799035	0.269202799838938\\
};
\addlegendentry{SFEM (tri subcell) (m=0.28)};

\end{axis}
\end{tikzpicture}%

%% file: Figures/platehole.pstex_t
\begin{picture}(0,0)%
\includegraphics{./Figures/platehole.pstex}%
\end{picture}%
\setlength{\unitlength}{3947sp}%
\begingroup\makeatletter\ifx\SetFigFont\undefined%
\gdef\SetFigFont#1#2#3#4#5{%
  \reset@font\fontsize{#1}{#2pt}%
  \fontfamily{#3}\fontseries{#4}\fontshape{#5}%
  \selectfont}%
\fi\endgroup%
\begin{picture}(4077,4240)(4051,-4874)
\put(8026,-2806){\makebox(0,0)[lb]{\smash{{\SetFigFont{14}{16.8}{\familydefault}{\mddefault}{\updefault}{\color[rgb]{0,0,0}$\sigma$}%
}}}}
\put(4246,-841){\makebox(0,0)[lb]{\smash{{\SetFigFont{14}{16.8}{\familydefault}{\mddefault}{\updefault}{\color[rgb]{0,0,0}$y$}%
}}}}
\put(4111,-4441){\makebox(0,0)[lb]{\smash{{\SetFigFont{12}{14.4}{\familydefault}{\mddefault}{\updefault}{\color[rgb]{0,0,0}$a$}%
}}}}
\put(7756,-4801){\makebox(0,0)[lb]{\smash{{\SetFigFont{14}{16.8}{\familydefault}{\mddefault}{\updefault}{\color[rgb]{0,0,0}$x$}%
}}}}
\put(5386,-2086){\makebox(0,0)[lb]{\smash{{\SetFigFont{12}{14.4}{\familydefault}{\mddefault}{\updefault}{\color[rgb]{0,0,0}$L$}%
}}}}
\put(6481,-3121){\rotatebox{90.0}{\makebox(0,0)[lb]{\smash{{\SetFigFont{12}{14.4}{\familydefault}{\mddefault}{\updefault}{\color[rgb]{0,0,0}$L$}%
}}}}}
\end{picture}%

%% file: Figures/pwHoledError.tikz
%
%
%
%
\begin{tikzpicture}

\begin{axis}[%
width=\figurewidth,
height=\figureheight,
scale only axis,
xmode=log,
xmin=100,
xmax=10000,
xminorticks=true,
xlabel={$\text{(total dof)}^{\text{1/2}}$},
ymode=log,
ymin=0.0001,
ymax=0.1,
yminorticks=true,
ylabel={Relative error in displacement},
legend style={draw=black,fill=white,legend cell align=left}
]
\addplot [
color=red,
solid,
mark=diamond,
mark options={solid}
]
table[row sep=crcr]{
404 0.007201360693307\\
790 0.004368235117637\\
1598 0.00186467533149\\
3178 0.0008014627252884\\
6360 0.000455433128974\\
};
\addlegendentry{Polygonal FEM (m=1.99)};

\addplot [
color=black,
solid,
mark=o,
mark options={solid}
]
table[row sep=crcr]{
404 0.021742978488148\\
790 0.012356149377212\\
1598 0.009242667275117\\
3178 0.00566477327497\\
6360 0.004641240777099\\
};
\addlegendentry{SFEM (tri subcell) (m=1.13)};

\addplot [
color=blue,
solid,
mark=square,
mark options={solid}
]
table[row sep=crcr]{
404 0.008762917977575\\
790 0.004484606720089\\
1598 0.00254281477256\\
3178 0.001354191955724\\
6360 0.000512273971356521\\
};
\addlegendentry{SFEM with stab (m=2.06)};

\end{axis}
\end{tikzpicture}%

%% file: Figures/LshapeError.tikz
%
%
%
%
\begin{tikzpicture}

\begin{axis}[%
width=\figurewidth,
height=\figureheight,
scale only axis,
xmode=log,
xmin=10,
xmax=1000,
xminorticks=true,
xlabel={$\text{(total dof)}^{1/2}$},
ymode=log,
ymin=0.01,
ymax=1,
yminorticks=true,
ylabel={Relative Error in displacement},
legend style={draw=black,fill=white,legend cell align=left}
]
\addplot [
color=blue,
solid,
mark=square,
mark options={solid}
]
table[row sep=crcr]{
17.8885438199983 0.138\\
34.8137903710584 0.0882\\
67.1416413263781 0.0761\\
137.564530312141 0.043\\
};
\addlegendentry{Polygonal FEM (m=0.51)};

\addplot [
color=red,
solid,
mark=o,
mark options={solid}
]
table[row sep=crcr]{
17.8885438199983 0.069507166475889\\
34.8137903710584 0.045136547806914\\
67.1416413263781 0.041743611828976\\
137.564530312141 0.026699084011734\\
};
\addlegendentry{SFEM with stab(m=0.47)};

\end{axis}
\end{tikzpicture}%

%% file: threeDexamples.tex
\subsubsection{Stability condition}
\begin{figure}[htpb]
\centering
\scalebox{0.65}{\input{./Figures/sfem3deig.pstex_t}}
\caption{Hexahedral element for the stability test}
\label{fig:hexelement}
\end{figure}
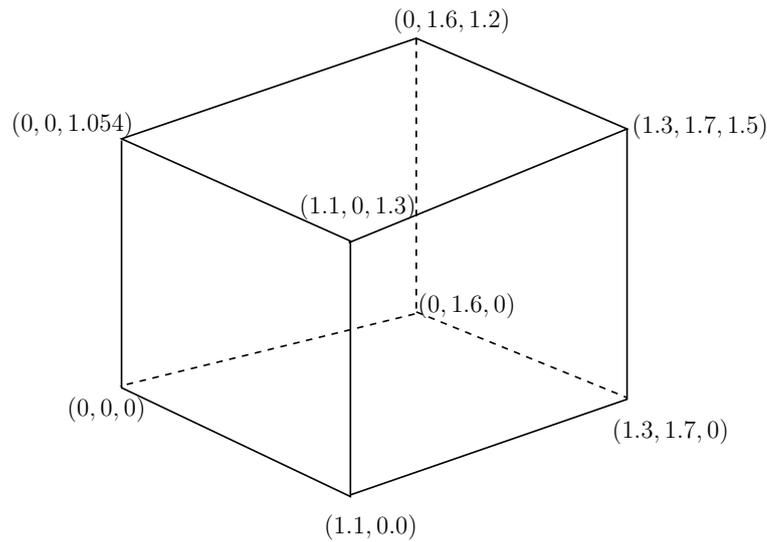

Before we proceed to study the accuracy and the convergence of the proposed method, we first investigate the stability condition by computing the eigenvalues of the hexahedral element shown in \fref{fig:hexelement}. The eigenvalues of the stiffness matrix computed by using the trilinear shape functions are $\lambda = \{0,0,0,0,0,0,1.47,3.11,3.55,4.28,4.82,6.04,6.44.6.98,8.75,9.85,\\
11.50,12.75,13.34,14.19,15.06,16.24,16.62,45.91 \}$. We solve the same problem with the proposed SFEM with and without stabilization. The eigenvalues of the stiffness matrix computed by SFEM without stabilization (i.e., with one subcell are): $\lambda= \{0,0,0,0,0,0,0,0,0,0,0,0,0,0,0,0,0,0,9.92,11.64,13.60,16.37,16.62,46.98\}$ and with stabilization ($\alpha=0.1$) are: $\lambda = \{ 0,0,0,0,0,0,9.32,9.45,10.42,10.64,11.30,11.49,11.51,\\ 11.51,11.51,11.51,11.51,11.51,11.53, 12.41,15.09,17.36,17.72,47.45\}$
It can be seen that SFEM with stabilization, like the FEM can capture the six zero energy modes corresponding to the physical rigid body modes. This indicates that the stiffness matrix is full rank and does not have any spurious energy modes. However, with one subcell and without stabilization, the eigenvalues of the stiffness matrix has 12 additional zero energy modes, these are non-physical. This is identical to the FEM with one integration point. Traditionally, for hexahedral elements, these are suppressed by stabilization procedures. In~\cite{bordasrabczuk2010,hunghiep2012}, the authors suppress the zero energy modes by adding the stiffness matrix computed with 8 subcells. In the present study, we add the stability term as given by \Eref{eqn:stabterm}.

\subsubsection{Patch Test}
\begin{figure}[htpb]
\centering
\includegraphics[scale=0.6]{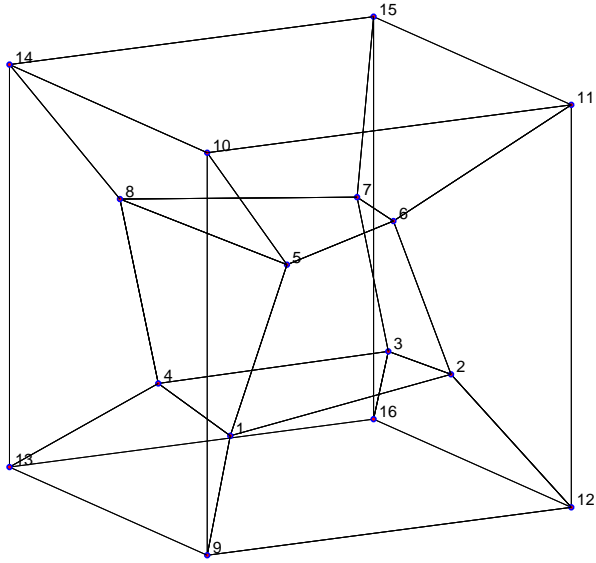}
\caption{Mesh used for displacement patch test: a cube containing warped elements (the coordinates of the interior nodes are given in Table \ref{tab:nodalcoord}).}
\label{fig:patchthreed}
\end{figure}

A three dimensional patch test with warped elements suggested in~\cite{macnealharder1985} is considered. The patch of elements shown in \fref{fig:patchthreed} is tested with the following displacement field applied on the outer boundary:
\begin{align}
u(x,y,z) &= \frac{1}{2}(2x + y + z)10^{-3} \nonumber \\
v(x,y,z) &= \frac{1}{2}(x +2y + z)10^{-3} \nonumber \\
w(x,y,z) &= \frac{1}{2}(x + y + 2z)10^{-3} 
\end{align}

The coordinates of the interior nodes are given in Table \ref{tab:nodalcoord}. The interior nodes are enclosed within a unit cube. It is observed that both the FEM with 3$\times$3$\times$3 Gaussian quadrature and cell-based smoothing technique with stabilization pass the patch test to machine precision.

\begin{table}[htpb]
\renewcommand{\arraystretch}{1}
\centering
\caption{Three dimensional patch test: coordinates of inner nodes. The inner nodes are enclosed within a unit cube~\cite{macnealharder1985}.}
\begin{tabular}{rrrr}
\hline
Node  & $x$ & $y$ & $z$ \\
number \\
\hline
1 & 0.249 & 0.342 & 0.192 \\
2 & 0.826 & 0.288 & 0.288 \\
3 & 0.850 & 0.649 & 0.263 \\
4 & 0.273 & 0.750 & 0.230 \\
5 & 0.320 & 0.186 & 0.643 \\
6 & 0.677 & 0.305 & 0.683 \\
7 & 0.788 & 0.693 & 0.644 \\
8 & 0.165 & 0.745 & 0.705 \\
\hline
\end{tabular}%
\label{tab:nodalcoord}%
\end{table}%

\subsubsection{Cantilever beam under shear load}
\begin{figure}[htpb]
\centering
\subfigure[50 elements]{\includegraphics[scale=0.2]{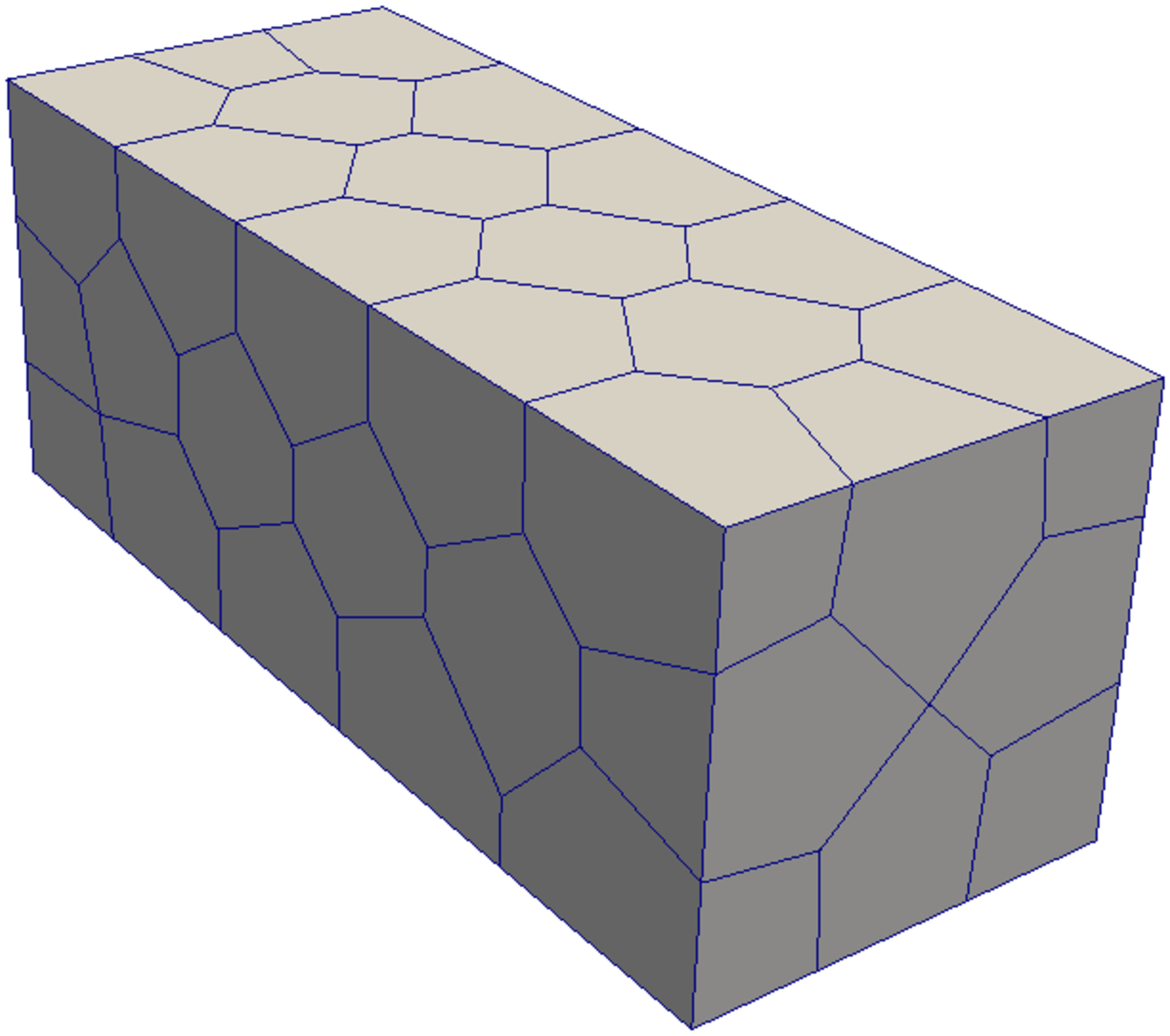}}
\subfigure[100 elements]{\includegraphics[scale=0.2]{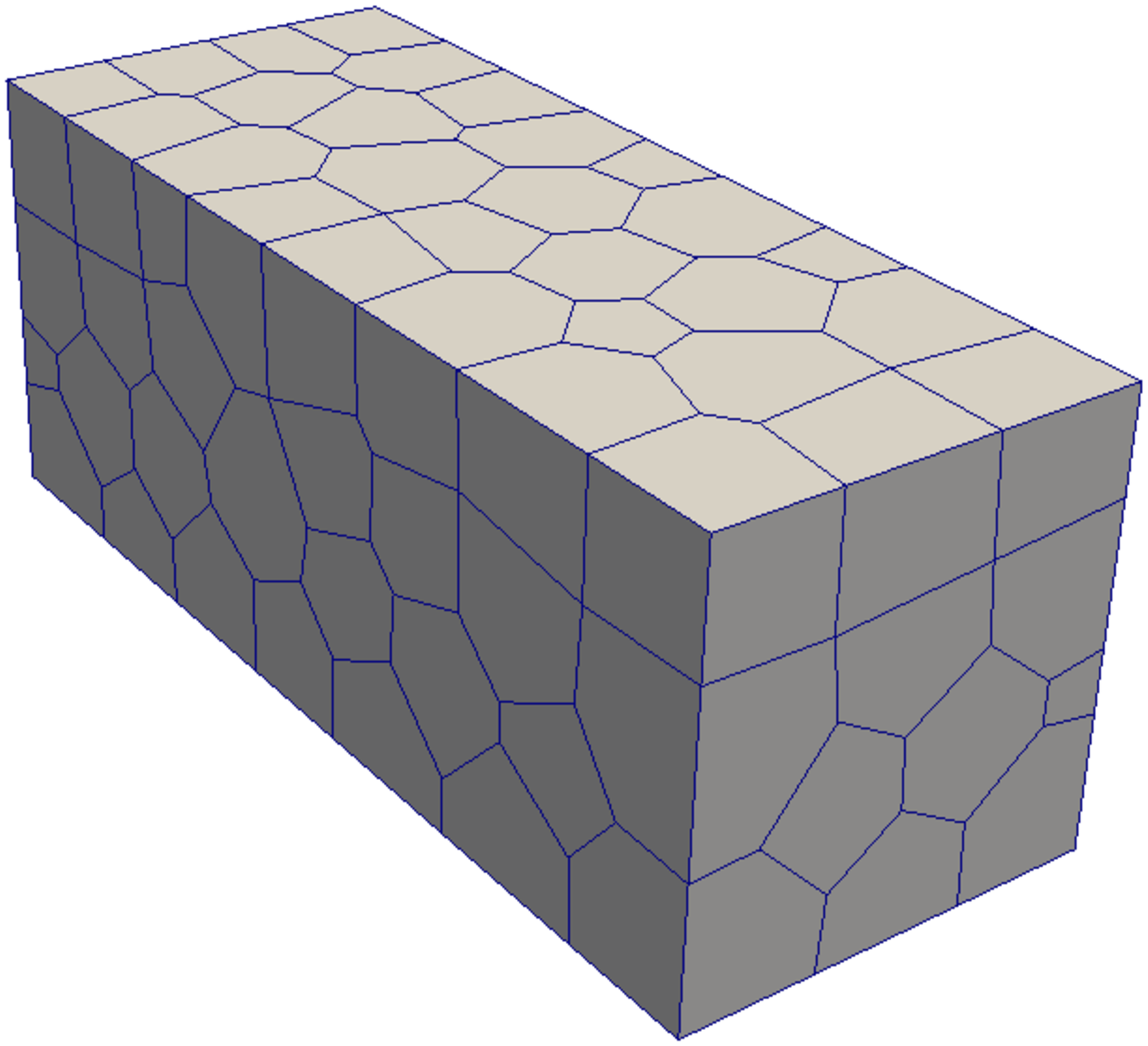}}
\subfigure[300 elements]{\includegraphics[scale=0.2]{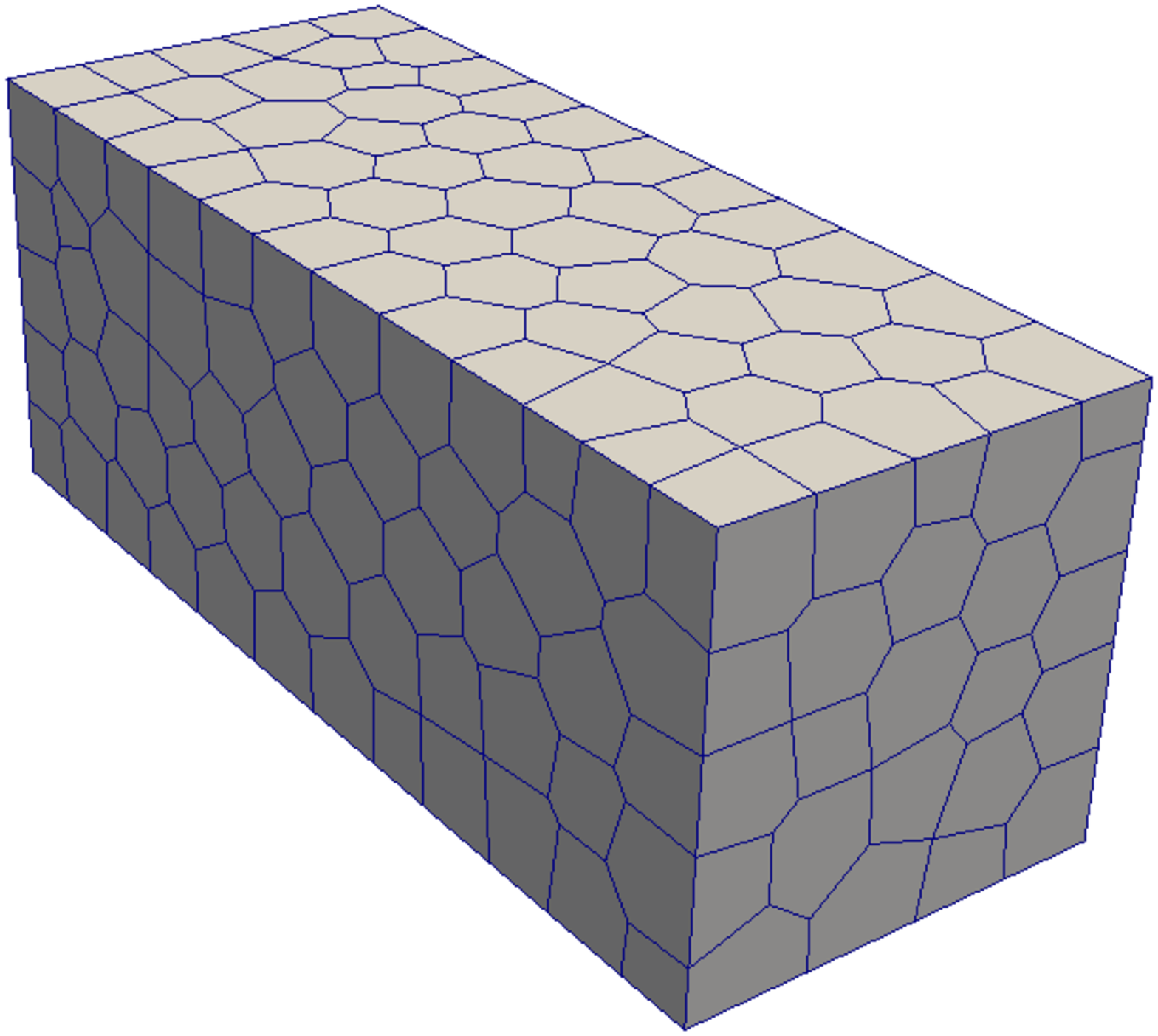}}
\subfigure[2000 elements]{\includegraphics[scale=0.2]{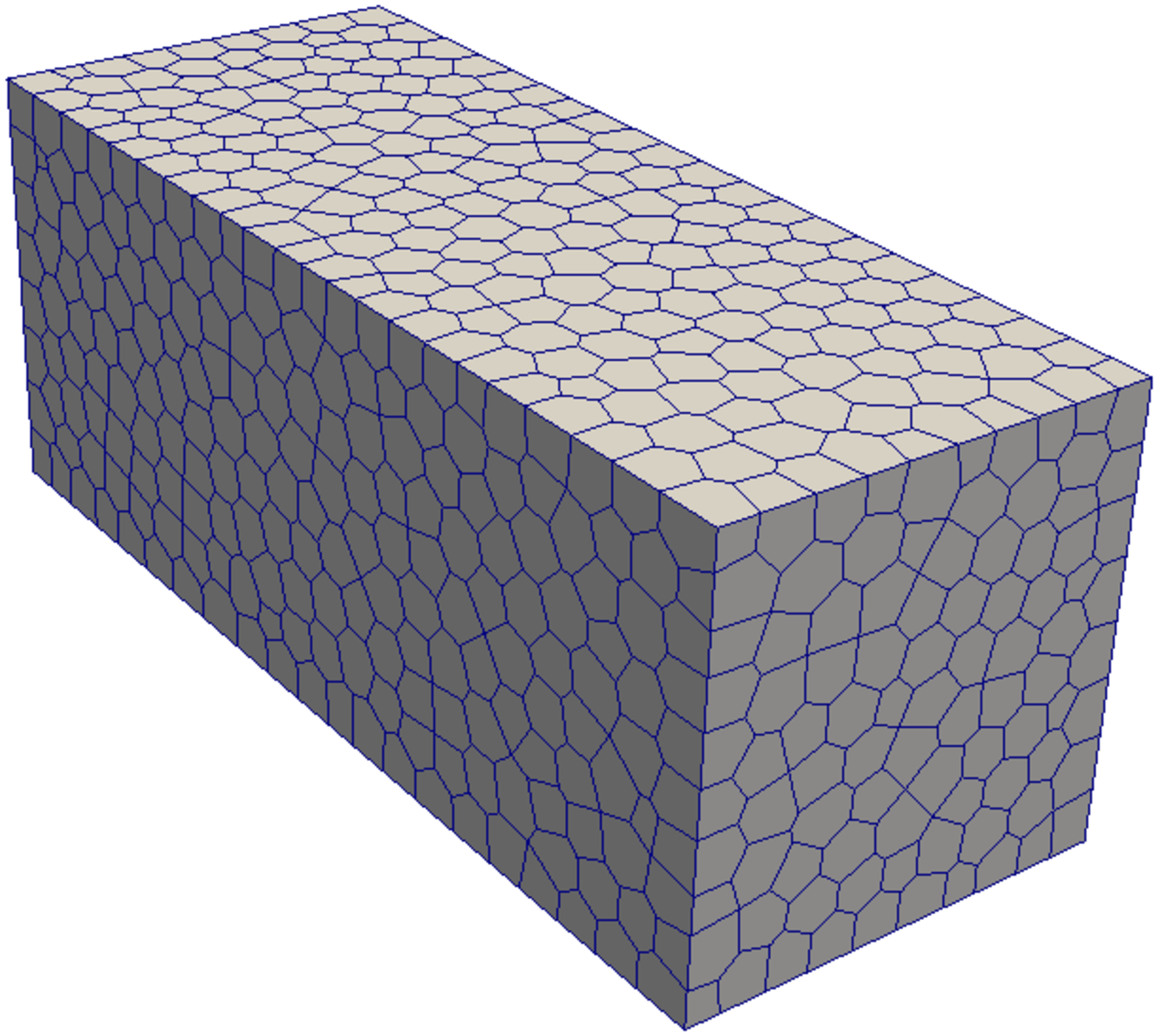}}
\caption{Random closed pack centroid Voronoi tessellation.}
\label{fig:cantirandmesh}
\end{figure}

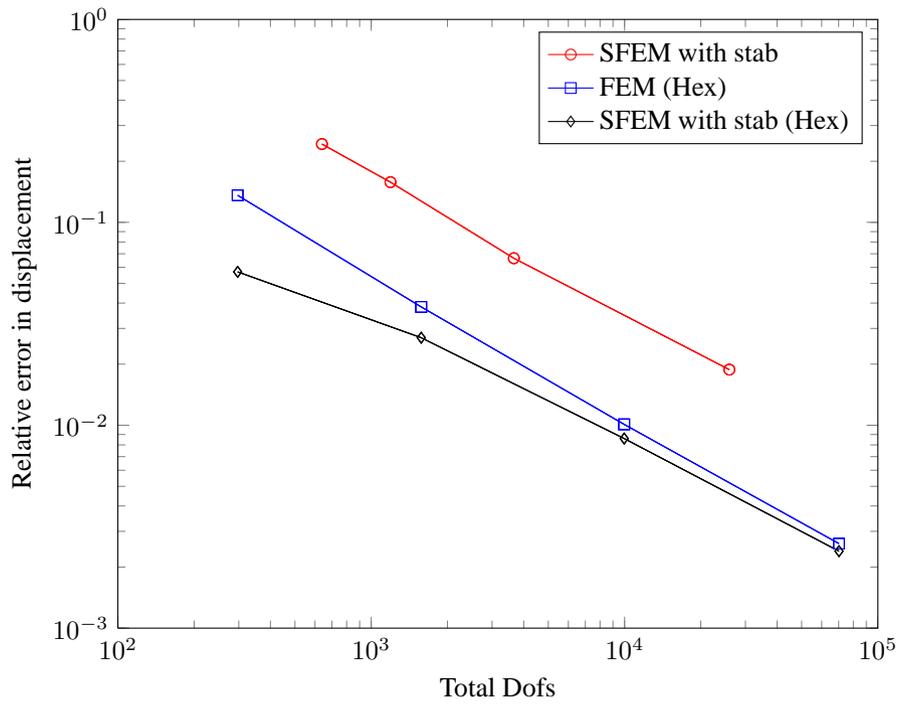
\begin{figure}[htpb]
\centering
\setlength\figureheight{8cm} 
\setlength\figurewidth{10cm}
\input{./Figures/cantibeam3DL2p.tikz}
\caption{Beam with shear load: convergence of the relative error in the displacement.}
\label{fig:cantibeam3dL2polygon}
\end{figure}

Consider a cantilever beam loaded in shear. The domain $\Omega$ for this problem is $[-1,1] \times [-1,1] \times [0,L]$. The material is assumed to be isotropic with Young's modulus, $E=$ 1 N/m$^2$ and Poisson's ratio $\nu=$ 0.3. The beam is subjected to a shear force $F$ at $z=L$ and at any cross section of the beam, the following conditions are satisfied:
\begin{equation}
\int\limits_{-a}^b \int\limits_{-a}^b  \sigma_{yz} ~\rmd x \rmd y = F \hspace{0.5cm} \int\limits_{-a}^b \int\limits_{-a}^b \sigma_{zz} y~\rmd x \rmd y = F z
\end{equation}
The Cauchy stress field is given by~\cite{barber2010}:
\begin{align}
\sigma_{xx}(x,y,z) &= \sigma_{xy}(x,y,z) = \sigma_{yy}(x,y,z)= 0 \nonumber \\
\sigma_{zz}(x,y,z) &= \frac{F}{I}yz \nonumber \\
\sigma_{xz}(x,y,z) &= \frac{2a^2 \nu F}{\pi^2 I (1+\nu)} \sum\limits_{n=0}^\infty \frac{ (-1)^n}{n^2} \sin \left(\frac{n \pi x}{a} \right) \frac{ \sinh \left( \frac{n \pi y}{a} \right)}{ \cosh \left( \frac{n \pi b}{a} \right)} \nonumber \\
\sigma_{yz}(x,y,z) &= \frac{(b^2-y^2)F}{2I} + \frac{ \nu F}{I(1+\nu)} \left[ \frac{ 3x^2-a^2}{6} - \frac{2a^2}{\pi^2} \sum \limits_{n=1}^\infty \frac{ (-1)^n}{n^2} \cos \left(\frac{n \pi x}{a} \right) \frac{ \cosh \left( \frac{n \pi y}{a} \right)}{ \cosh \left( \frac{n \pi b}{a} \right)} \right]
\end{align}
The corresponding displacement field is given by~\cite{bishop2014}:
\begin{align}
u(x,y,z) &= - \frac{\nu F}{EI} xyz \nonumber \\
v(xy,z) &= \frac{F}{EI} \left[ \frac{\nu (x^2-y^2)z}{2} - \frac{z^3}{6} \right] \nonumber \\
w(x,y,z) &= \frac{F}{EI} \left[ \frac{y (\nu x^2+z^2)}{2} + \frac{\nu y^3}{6} + (1+\nu) \left( b^2y-\frac{y^3}{3} \right) - \frac{\nu a^2 y}{3} - \frac{4 \nu a^3}{\pi^3} \sum\limits_{n=0}^\infty \frac{ (-1)^n}{n^2} \cos \left(\frac{n \pi x}{a} \right) \frac{ \sinh \left( \frac{n \pi y}{a} \right)}{ \cosh \left( \frac{n \pi b}{a} \right)} \right]
\label{eqn:3dcantianadisp}
\end{align}
where $E$ is the Young's modulus, $\nu$ is Poisson's ratio and $I= 4ab^3/3$ is the second moment of area about the $x$-axis. Two types of meshes are considered: (a) a regular hexahedral mesh and (2) a random closed -pack Voronoi mesh. Four levels of mesh refinement are considered for both hexahedral meshes (2$\times$2$\times$10, 4$\times$4$\times$20, 8$\times$8$\times$40, 16$\times$16$\times$80) and for the random Voronoi meshes. \fref{fig:cantirandmesh} shows the random Voronoi mesh employed for this study. The length of the beam is $L=$ 5 and the shear load is taken as $F=$ 1. Analytical displacements given by \Eref{eqn:3dcantianadisp} are applied at $z=L$, whilst the beam is loaded in shear at $z=0$. All other faces are assumed to be traction free. \fref{fig:cantibeam3dL2polygon} shows the relative error in the displacement norm with mesh refinement. It can be seen that both formulations, viz., FEM and SFEM with stabilization converge to the analytical solution with mesh refinement. It is also observed that the SFEM with stabilization yields slightly more accurate results when compared to the FEM with full integration. Also, shown in \fref{fig:cantibeam3dL2polygon} is the convergence of the method with polyhedra meshes. Although, the formulation when applied to polyhedra elements, converges with mesh refinement, it is not as accurate as the hexahedral elements. This can attributed to the parameter $\alpha^\ast$ employed in this study and this observation is consistent with the results reported in the literature~\cite{gaintalischi2014}.

\begin{figure}[htpb]
\centering
\subfigure[]{\scalebox{0.55}{\input{./Figures/threedbeam.pstex_t}}}
\subfigure[]{\includegraphics[scale=0.30]{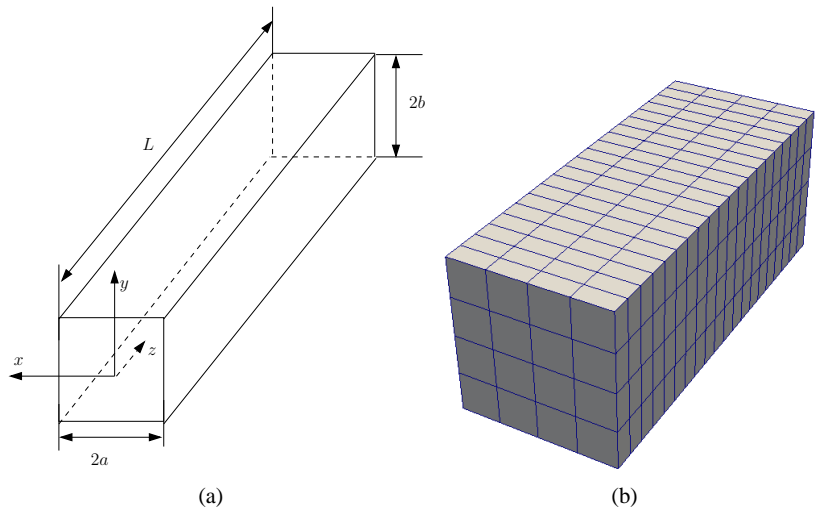}}
\caption{Cantilever beam: (a) Geometry, length $L$ and rectangular cross-section of wdith $2a$ and height $2b$. For the present study, the following dimensions are considered: $L=$ 5, $a=b=$ 1 and (b) A structured hexahedral mesh (4$\times$4$\times$20.}
\label{fig:canti3dbeam}
\end{figure}

%% file: Figures/sfem3deig.pstex_t
\begin{picture}(0,0)%
\includegraphics{./Figures/sfem3deig.pstex}%
\end{picture}%
\setlength{\unitlength}{4144sp}%
\begingroup\makeatletter\ifx\SetFigFont\undefined%
\gdef\SetFigFont#1#2#3#4#5{%
  \reset@font\fontsize{#1}{#2pt}%
  \fontfamily{#3}\fontseries{#4}\fontshape{#5}%
  \selectfont}%
\fi\endgroup%
\begin{picture}(5685,4882)(2146,-5993)
\put(2161,-2251){\makebox(0,0)[lb]{\smash{{\SetFigFont{14}{16.8}{\familydefault}{\mddefault}{\updefault}{\color[rgb]{0,0,0}$(0,0,1.054)$}%
}}}}
\put(2671,-4846){\makebox(0,0)[lb]{\smash{{\SetFigFont{14}{16.8}{\familydefault}{\mddefault}{\updefault}{\color[rgb]{0,0,0}$(0,0,0)$}%
}}}}
\put(5011,-5911){\makebox(0,0)[lb]{\smash{{\SetFigFont{14}{16.8}{\familydefault}{\mddefault}{\updefault}{\color[rgb]{0,0,0}$(1.1,0.0)$}%
}}}}
\put(7636,-5041){\makebox(0,0)[lb]{\smash{{\SetFigFont{14}{16.8}{\familydefault}{\mddefault}{\updefault}{\color[rgb]{0,0,0}$(1.3,1.7,0)$}%
}}}}
\put(7816,-2266){\makebox(0,0)[lb]{\smash{{\SetFigFont{14}{16.8}{\familydefault}{\mddefault}{\updefault}{\color[rgb]{0,0,0}$(1.3,1.7,1.5)$}%
}}}}
\put(5641,-1306){\makebox(0,0)[lb]{\smash{{\SetFigFont{14}{16.8}{\familydefault}{\mddefault}{\updefault}{\color[rgb]{0,0,0}$(0,1.6,1.2)$}%
}}}}
\put(5866,-3901){\makebox(0,0)[lb]{\smash{{\SetFigFont{14}{16.8}{\familydefault}{\mddefault}{\updefault}{\color[rgb]{0,0,0}$(0,1.6,0)$}%
}}}}
\put(4786,-3001){\makebox(0,0)[lb]{\smash{{\SetFigFont{14}{16.8}{\familydefault}{\mddefault}{\updefault}{\color[rgb]{0,0,0}$(1.1,0,1.3)$}%
}}}}
\end{picture}%

%% file: Figures/threedbeam.pstex_t
\begin{picture}(0,0)%
\includegraphics{./Figures/threedbeam.pstex}%
\end{picture}%
\setlength{\unitlength}{3947sp}%
\begingroup\makeatletter\ifx\SetFigFont\undefined%
\gdef\SetFigFont#1#2#3#4#5{%
  \reset@font\fontsize{#1}{#2pt}%
  \fontfamily{#3}\fontseries{#4}\fontshape{#5}%
  \selectfont}%
\fi\endgroup%
\begin{picture}(4704,5256)(3049,-5095)
\put(4636,-3781){\makebox(0,0)[lb]{\smash{{\SetFigFont{12}{14.4}{\familydefault}{\mddefault}{\updefault}{\color[rgb]{0,0,0}$z$}%
}}}}
\put(7591,-991){\makebox(0,0)[lb]{\smash{{\SetFigFont{12}{14.4}{\familydefault}{\mddefault}{\updefault}{\color[rgb]{0,0,0}$2b$}%
}}}}
\put(3991,-5026){\makebox(0,0)[lb]{\smash{{\SetFigFont{12}{14.4}{\familydefault}{\mddefault}{\updefault}{\color[rgb]{0,0,0}$2a$}%
}}}}
\put(4321,-3031){\makebox(0,0)[lb]{\smash{{\SetFigFont{12}{14.4}{\familydefault}{\mddefault}{\updefault}{\color[rgb]{0,0,0}$y$}%
}}}}
\put(3121,-3901){\makebox(0,0)[lb]{\smash{{\SetFigFont{12}{14.4}{\familydefault}{\mddefault}{\updefault}{\color[rgb]{0,0,0}$x$}%
}}}}
\put(4576,-1471){\makebox(0,0)[lb]{\smash{{\SetFigFont{12}{14.4}{\familydefault}{\mddefault}{\updefault}{\color[rgb]{0,0,0}$L$}%
}}}}
\end{picture}%

%% file: Figures/cantibeam3DL2p.tikz
%
%
%
%
\begin{tikzpicture}

\begin{axis}[%
width=\figurewidth,
height=\figureheight,
scale only axis,
xmode=log,
xmin=100,
xmax=100000,
xminorticks=true,
xlabel={Total Dofs},
ymode=log,
ymin=0.001,
ymax=1,
yminorticks=true,
ylabel={Relative error in displacement},
legend style={draw=black,fill=white,legend cell align=left}
]
\addplot [
color=red,
solid,
mark=o,
mark options={solid}
]
table[row sep=crcr]{
639 0.24313157971727\\
1191 0.157760842819165\\
3657 0.066543426538196\\
25929 0.018794375418417\\
};
\addlegendentry{SFEM with stab};

\addplot [
color=blue,
solid,
mark=square,
mark options={solid}
]
table[row sep=crcr]{
297 0.1359\\
1575 0.0383\\
9963 0.0101\\
70227 0.002609469\\
};
\addlegendentry{FEM (Hex)};

\addplot [
color=black,
solid,
mark=diamond,
mark options={solid}
]
table[row sep=crcr]{
297 0.057\\
1575 0.027\\
9963 0.0086\\
70227 0.002394743\\
};
\addlegendentry{SFEM with stab (Hex)};

\addplot [
color=red,
solid,
mark=o,
mark options={solid},
forget plot
]
table[row sep=crcr]{
639 0.24313157971727\\
1191 0.157760842819165\\
3657 0.066543426538196\\
25929 0.018794375418417\\
};
\addplot [
color=blue,
solid,
mark=square,
mark options={solid},
forget plot
]
table[row sep=crcr]{
297 0.1359\\
1575 0.0383\\
9963 0.0101\\
70227 0.002609469\\
};
\addplot [
color=black,
solid,
mark=diamond,
mark options={solid},
forget plot
]
table[row sep=crcr]{
297 0.057\\
1575 0.027\\
9963 0.0086\\
70227 0.002394743\\
};
\addplot [
color=red,
solid,
mark=o,
mark options={solid},
forget plot
]
table[row sep=crcr]{
639 0.24313157971727\\
1191 0.157760842819165\\
3657 0.066543426538196\\
25929 0.018794375418417\\
};
\addplot [
color=blue,
solid,
mark=square,
mark options={solid},
forget plot
]
table[row sep=crcr]{
297 0.1359\\
1575 0.0383\\
9963 0.0101\\
70227 0.002609469\\
};
\addplot [
color=black,
solid,
mark=diamond,
mark options={solid},
forget plot
]
table[row sep=crcr]{
297 0.057\\
1575 0.027\\
9963 0.0086\\
70227 0.002394743\\
};
\end{axis}
\end{tikzpicture}%

%% file: ApplicationLefm.tex
\input{sbfem_crk}

\subsubsection{Plate with double edge crack in tension}
\begin{figure}[htpb]
\centering
\subfigure[]{\scalebox{0.65}{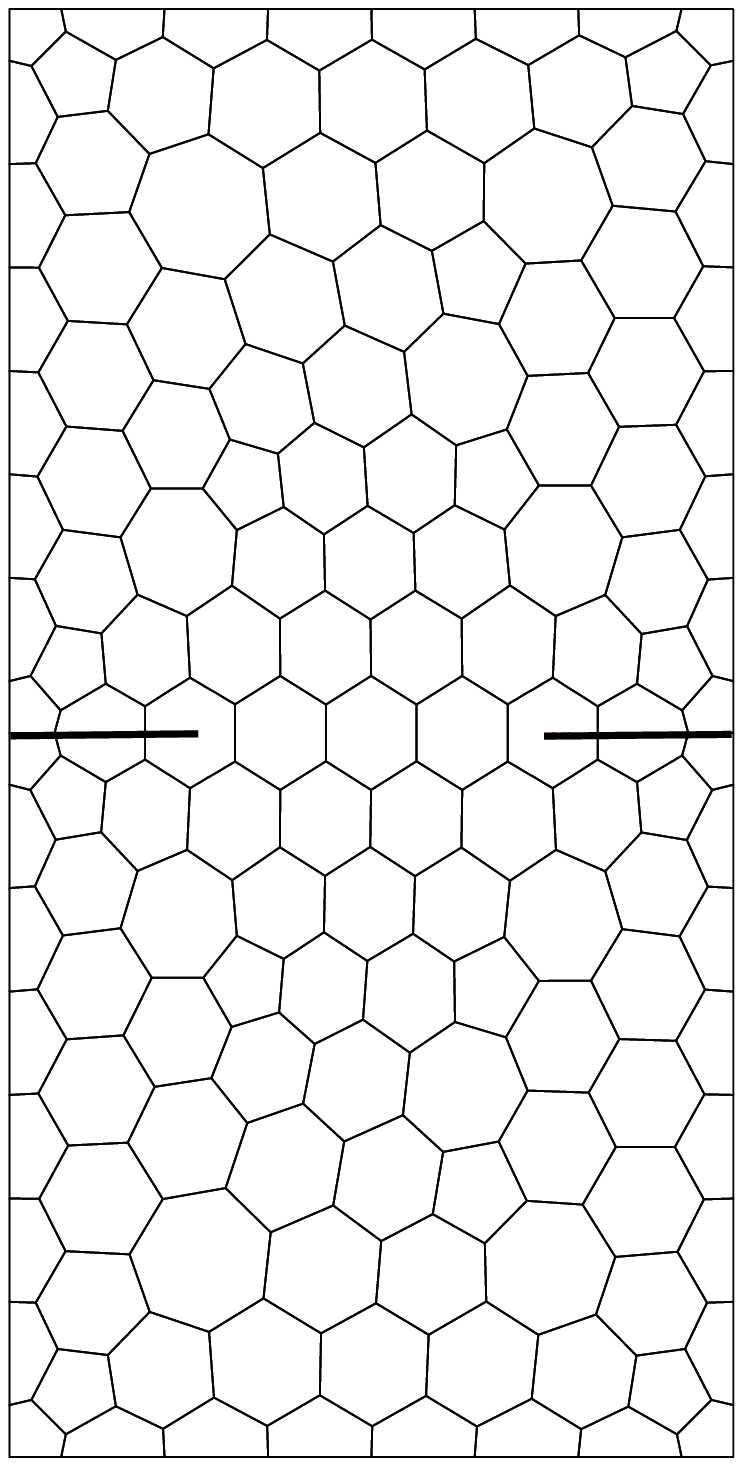}}
\subfigure[]{\includegraphics[scale=0.45]{./Figures/dedgcrk}}
\caption{Plate with an edge under tension: (a) geometry and boundary conditions and (b) domain discretized with polygonal elements.}
\label{fig:edgcrkgeom}
\end{figure}

The plate with double edge crack subjected to a uniform tension at both ends as shown in \fref{fig:edgcrkgeom} is considered. In the computations, the ratio of the crack length, $a$, to the width of the plate, $H$, is $a/H=$ 0.25. The material properties of the plate are: Young's modulus, $E=$ 200 GPa and Poisson's ratio, $\nu=0.3$. In this example, plane stress conditions are assumed. The empirical mode I SIF that is given by:
\begin{equation}
K_{\rm{I}}^{\mathrm{ref}}=C\sigma\sqrt{\pi a}\label{eq: KI ref}
\end{equation}
where $C$ is a correction factor. For $a/b > 0.4, b = H/2$, the correction factor is given by~\cite{tadaparis2000}:
\begin{equation}
C=1.12+0.203\left(\frac{a}{b}\right)-1.197\left(\frac{a}{b}\right)^{2}+1.930\left(\frac{a}{b}\right)^{3} 
\label{eqn:correctionfactor1}
\end{equation}

The above factor corrects for an infinite plate with an accuracy of 2$\%$. For the chosen parameters, the reference normalized SIF is $K_I/\sqrt{\pi a}=$ 1.1635. The plate is discretized with a polygonal mesh. For the polygons containing the crack tip, we employ the SBFEM technique to capture the singularity. In this polygon, each edge is further discretized with 5 linear elements so that the angular variation of the SIF can be computed accurately~\cite{ooisong2012}. For the elements that do not contain the crack tip, we employ the SFEM with stabilization to compute the stiffness matrix. The convergence of the mode I SIF with mesh refinement is give in Table \ref{tab:modeIedgeCrk}. It can be seen that the proposed method converges to the empirical relation with mesh refinement.

\begin{table}[htpb]
\centering
\caption{Plate with an edge crack in remote tension: convergence of mode I SIF.}
\begin{tabular}{rrrrr}
\hline
$h$ & Number &  Number & $K_I$ &  $K_I/\sqrt{\pi a}$\\
& of Polygons & of nodes & \\
\hline
0.25 & 32 & 111 &  1.0524 &  1.1875 \\
0.125 & 91 & 247  & 1.0555  & 1.1910 \\
0.0625 & 332 & 777  & 1.0452 &  1.1794 \\
0.03125 & 1229 & 2659  & 1.0431 &  1.1770 \\
0.015625 & 4940 & 10245  & 1.0432  & 1.1772 \\
\hline
\end{tabular}
\label{tab:modeIedgeCrk}
\end{table}

\subsubsection{Angled crack in an isotropic material}
In this example, a plate with an angled crack subjected to far field bi-axial stress field, $\bvsig$ (see \fref{fig:inclcrkgeom}) with $a/w=$ 0.1, $\sigma_1=$ 1 and $\sigma_2=$ 2 is considered. In this example, the mode I and the mode II SIFs, $K_{\rm I}$ and $K_{\rm II}$, respectively, are obtained as a function of the crack angle $\beta$. For the loads shown, the analytical SIF for an infinite plate are given by~\cite{aliabadirooke1987}:
\begin{align}
K_{\rm I} &= (\sigma_2 \sin^2 \beta + \sigma_1 \cos^2 \beta) \sqrt{\pi a} \nonumber \\
K_{\rm II} &= (\sigma_2 - \sigma_1) \sin \beta \cos \beta \sqrt{\pi a}
\label{eqn:inclcrcempirical}
\end{align}
The material properties of the plate are: Young's modulus, $E=$ 200 GPa and Poisson's ratio, $\nu=0.3$.  In this example, the plate is discretized with polygon meshes (containing 300 polygons). The variations of the mode I and mode II SIF with the crack orientation $\beta$ are presented in Table \ref{table:inclinedcracksifs}. It can be observed that the results from the present method agree very well with the reference solution.
\begin{figure}[htpb]
\centering
\scalebox{0.7}{\input{./Figures/inclCrk.pstex_t}}
\caption{Plate with an oblique crack: geometry and boundary conditions.}
\label{fig:inclcrkgeom}
\end{figure}

\begin{table}[htpb]
\centering
\renewcommand{\arraystretch}{1.5}
\caption{Mode I and Mode II SIF for a plate with an inclined crack.}
\begin{tabular}{lrrrrrrr}
\hline 
$\beta$ & \multicolumn{3}{c}{mode I SIF} && \multicolumn{3}{c}{mode II SIF} \\
\cline{2-4}\cline{6-8}
 & \Eref{eqn:inclcrcempirical} & Crack Tip A & Crack Tip B && \Eref{eqn:inclcrcempirical} & Crack Tip A & Crack Tip B \\
\hline
0$^\circ$ & 1.0000 & 1.0176 & 1.0168 && 0.0000 & 0.0000 & 0.0000 \\
15$^\circ$ & 1.0670 & 1.0937 & 1.0876 && 0.2500 & 0.2343 & 0.2453 \\
30$^\circ$ & 1.2500 & 1.2786 & 1.2786 && 0.4330 & 0.4380 & 0.4379 \\
45$^\circ$ & 1.5000 & 1.5281 & 1.5266 && 0.5000 & 0.5039 & 0.5053 \\
60$^\circ$ & 1.7500 & 1.7893 & 1.7893 && 0.4330 & 0.4427 & 0.4429 \\
75$^\circ$ & 1.9330 & 1.9855 & 1.9738 && 0.2500 & 0.2880 & 0.2669 \\
90$^\circ$ & 2.0000 & 2.0351 & 2.0336 && 0.0000 & 0.0018 & 0.0122 \\
\hline
\end{tabular}
\label{table:inclinedcracksifs}
\end{table}


%% file: sbfem_crk.tex
\subsection{Application to linear elastic fracture mechanics}
\label{sfemsbfem}
The SFEM with stabilization discussed above over arbitrary polygons and polyhedra can be applied to problems with strong discontinuity and singularity. However, to accurately capture the asymptotic fields at the crack tip, a very fine mesh in combination with singular elements at the crack tip is usually required. This poses additional difficulties when the crack evolves. Another possibility is to enrich the approximation space with functions that can capture the discontinuity and singularity~\cite{melenkbabuvska1996,belytschkogracie2009}. In the literature, the latter method is referred to as the Generalized FEM (GFEM)/extended finite element method (XFEM). In~\cite{bordasnatarajan2011}, the authors combined the strain smoothing with the XFEM. It was observed that in the case of enrichment schemes for linear elastic fracture mechanics, the method yields less accurate results compared to the conventional XFEM. However, for the elements that are completely intersected by the discontinuous surface, with the strain smoothing operation, further sub-division is not required. 

\begin{figure}[htpb]
\centering
\includegraphics[scale=0.75]{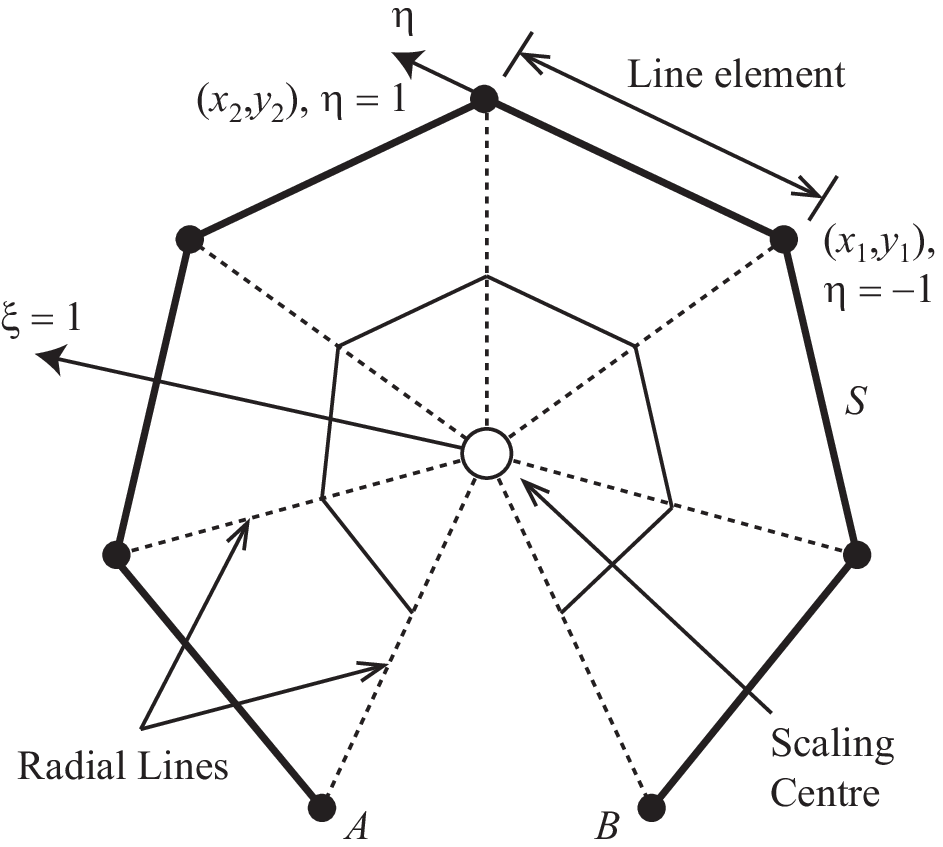}
\caption{Cracked polygon representation by the scaled boundary finite element method.}
\label{fig:Polygon-representation-by}
\end{figure}

In this study, we propose to couple the SFEM with the scaled boundary finite element method (SBFEM) to model problems with strong discontinuity and singularities. The SBFEM is a novel method that has the advantages of both the FEM and the boundary element method (BEM). As in the FEM, no fundamental solution is required and as in BEM, the problem dimension is reduced by one. The SBFEM is a semi-analytical method and relies on defining a `scaling centre' from which the entire boundary is visible. This is similar to the concept of `star convexity'. The boundary is divided into conventional linear finite elements, whilst the solution is sought analytically in the radial direction~\cite{wolfsong2000}.  Moreover, by exploiting the special characteristics of the scaling centre,  the stress intensity factors can be computed directly. When modelling a crack/notched surface the scaling centre is placed at the crack tip. The straight crack/notch edges are formed by scaling the nodes $A$ and $B$ on the boundary and the crack surfaces are not discretized (see ~\fref{fig:Polygon-representation-by}).

\paragraph{Displacement approximation} The geometry of the element described by the coordinates on the boundary $\xx_b(\eta)$ is expressed as:
\begin{equation}
\xx_b(\eta) = \mathbf{N}(\eta) \xx_b
\end{equation}
where $\mathbf{N}(\eta)$ is the shape function matrix of the finite elements discretising the polygon boundary. The standard 1D Gauss-Lobatto shape functions or Lagrange shape functions can be used. In this study, we employ Lagrange shape functions. The displacements of a point in a polygon is approximated by: 
\begin{equation}
\mathbf{u}(\xi,\eta)=\mathbf{N}(\eta)\mathbf{u}(\xi)\label{eqn:dispapprox}
\end{equation}
where $\mathbf{u}(\xi)$ are radial displacement functions. Substituting \Eref{eqn:dispapprox} in the definition of strain-displacement relations, the strains $\boldsymbol{\varepsilon}(\xi,\eta)$ are expressed as:
\begin{equation}
\boldsymbol{\varepsilon}(\xi,\eta)=\mathbf{L}\mathbf{u}(\xi,\eta)\label{eqn:sbfemstrain}
\end{equation}
 where $\mathbf{L}$ is a linear operator matrix formulated in the
scaled boundary coordinates as 
\begin{equation}
\mathbf{L}=\mathbf{b}_{1}(\eta)\frac{\partial}{\partial\xi}+\xi^{-1}\mathbf{b}_{2}(\eta)\label{eqn:Loperator}
\end{equation}
 with 
\begin{align}
\mathbf{b}_{1}(\eta) & =\frac{1}{|\mathbf{J}(\eta)|}\left[\begin{array}{cc}
y_{\eta}(\eta)_{,\eta} & 0\\
0 & -x_{\eta}(\eta)_{,\eta}\\
-x_{\eta}(\eta)_{,\eta} & y_{\eta}(\eta)_{,\eta}
\end{array}\right] \nonumber \\
\mathbf{b}_{2}(\eta) & =\frac{1}{|\mathbf{J}(\eta)|}\left[\begin{array}{cc}
-y_{\eta}(\eta) & 0\\
0 & x_{\eta}(\eta)\\
x_{\eta}(\eta) & y_{\eta}(\eta)
\end{array}\right]\label{eq:b2}
\end{align}
By following the procedure outlined in~\cite{wolfsong2001,deekswolf2002}, the following ODE is obtained:
\begin{equation}
\mathbf{E}_{0}\xi^{2}\mathbf{u}(\xi)_{,\xi\xi}+(\mathbf{E}_{0}+\mathbf{E}_{1}^{\mathrm{T}}-\mathbf{E}_{1})\xi\mathbf{u}(\xi)_{,\xi}-\mathbf{E}_{2}\mathbf{u}(\xi)=0\label{eqn:governODEsbfem}
\end{equation}
where $\mathbf{E}_{0},\mathbf{E}_{1}$ and $\mathbf{E}_{2}$ are coefficient matrices given by:
\begin{align}
\mathbf{E}_{0} & =\int_{\eta}\mathbf{B}_{1}(\eta)^{{\rm T}}\mathbf{D}\mathbf{B}_{1}(\eta)|\mathbf{J}(\eta)|d\eta,\nonumber \\
\mathbf{E}_{1} & =\int_{\eta}\mathbf{B}_{2}(\eta)^{{\rm T}}\mathbf{D}\mathbf{B}_{1}(\eta)|\mathbf{J}(\eta)|d\eta,\nonumber \\
\mathbf{E}_{2} & =\int_{\eta}\mathbf{B}_{2}(\eta)^{{\rm T}}\mathbf{D}\mathbf{B}_{2}(\eta)|\mathbf{J}(\eta)|d\eta.\label{eqn:coeffmat}
\end{align}
Using ~\Eref{eqn:sbfemstrain} and Hooke's law $\boldsymbol{\sigma}=\mathbf{D}\boldsymbol{\varepsilon}$, the stresses $\boldsymbol{\sigma}(\xi,\eta)$ is expressed as 
\begin{equation}
\boldsymbol{\sigma}(\xi,\eta)=\mathbf{D}\left(\mathbf{B}_{1}(\eta)\mathbf{u}(\xi)_{,\xi}+\xi^{-1}\mathbf{B}_{2}(\eta)\mathbf{u}(\xi)\right)\label{eqn:sbfemstress}
\end{equation}
where $\mathbf{D}$ is the material constitutive matrix and the determinant of the Jacobian matrix is:
\begin{equation}
|\mathbf{J}(\eta)| = x_b(\eta) y_b(\eta)_{,\eta} - y_b(\eta)x_b(\eta)_{,\eta}
\end{equation}
The coefficient matrices are evaluated element-by-element on the polygon boundary and assembled over a polygon. This process is similar to the standard FE procedure of assemblage. \Eref{eqn:governODEsbfem} is a homogeneous second-order ordinary differential equation. Its solution is obtained by introducing the variable $\boldsymbol{\chi}(\xi)$
\begin{align}
\boldsymbol{\chi}(\xi)= & \left\{ \begin{array}{c}
\mathbf{u}(\xi)\\
\mathbf{q}(\xi)
\end{array}\right\} \label{eq:chi}
\end{align}
 where $\mathbf{q}(\xi)$ is the internal load vector 
\begin{align}
\mathbf{q}(\xi)= & \mathbf{E}_{0}\xi\mathbf{u}(\xi)_{,\xi}+\mathbf{E}_{1}^{{\rm T}}\mathbf{u}(\xi)
\end{align}
The boundary nodal forces are related to the displacement functions by: 
\begin{equation}
\mathbf{f}=\mathbf{q}(\xi=1)=(\mathbf{E}_{0}\xi\mathbf{u}(\xi)_{,\xi}+\mathbf{E}_{1}^{{\rm T}}\mathbf{u}(\xi))|_{\xi=1}\label{eqn:nodalforce}
\end{equation}
This allows~\Eref{eqn:governODEsbfem} to be transformed into a first order ordinary differential equation with twice the number of unknowns in an element as: 
\begin{equation}
\xi\boldsymbol{\chi}(\xi)_{,\xi}=-\mathbf{Z}\boldsymbol{\chi}(\xi)\label{eq:first order}
\end{equation}
 where $\mathbf{Z}$ is a Hamiltonian matrix 
\begin{equation}
\mathbf{Z}=\left[\begin{array}{cc}
\mathbf{E}_{0}^{-1}\mathbf{E}_{1}^{\mathrm{T}} & -\mathbf{E}_{0}^{-1}\\
\mathbf{E}_{1}\mathbf{E}_{0}^{-1}\mathbf{E}_{1}^{\mathrm{T}}-\mathbf{E}_{2} & -\mathbf{E}_{1}\mathbf{E}_{0}^{-1}
\end{array}\right]\label{eq:Hamiltonian matrix}
\end{equation}
An eigenvalue decomposition of $\mathbf{Z}$ is performed. The blocks of eigenvalues and transformation matrices necessary are:
\begin{align}
\mathbf{Z}\left[\begin{array}{c}
\boldsymbol{\Phi}_{\mathrm{u}}\\
\boldsymbol{\Phi}_{\mathrm{q}}
\end{array}\right] & =\left[\begin{array}{c}
\boldsymbol{\Phi}_{\mathrm{u}}\\
\boldsymbol{\Phi}_{\mathrm{q}}
\end{array}\right]\boldsymbol{\Lambda}_{\mathrm{n}}\label{eq:eigen decomp}
\end{align}
In \Eref{eq:eigen decomp}, $\boldsymbol{\Lambda}_{\mathrm{n}}=\mathrm{diag}\left(\lambda_{1},\,\lambda_{2},\,...,\lambda_{n}\right)$ contains the eigenvalues with negative real part. $\boldsymbol{\Phi}_{\mathrm{u}}$ and $\boldsymbol{\Phi}_{\mathrm{q}}$ are the corresponding transformation matrices of $\boldsymbol{\Lambda}_{\mathrm{n}}$. They represent the modal displacements and forces, respectively. The general solution of \Eref{eq:first order} is given by:
\begin{align}
\mathbf{u}(\xi)= & \boldsymbol{\Phi}_{\mathrm{u}}\xi^{-\boldsymbol{\Lambda}_{\mathrm{n}}}\mathbf{c}\\
\mathbf{q}(\xi)= & \boldsymbol{\Phi}_{\mathrm{q}}\xi^{-\boldsymbol{\Lambda}_{\mathrm{n}}}\mathbf{c}
\label{eqn:raddispfunc}
\end{align}
where $\mathbf{c}$ are integration constants that are obtained from the nodal displacements $\mathbf{u}_{\mathrm{b}}=\mathbf{u}(\xi=1)$ as:
\begin{align}
\mathbf{c}= & \boldsymbol{\Phi}_{\mathrm{u}}^{-1}\mathbf{u}_{\mathrm{b}}\label{eq:int constants}
\end{align}
The complete displacement field of a point defined by the sector covered by a line element on the element is obtained by substituting~\Eref{eqn:raddispfunc} into \Eref{eqn:dispapprox} resulting in:
\begin{align}
\mathbf{u}(\xi,\eta)= & \mathbf{N}(\eta)\boldsymbol{\Phi}_{\mathrm{u}}\xi^{-\boldsymbol{\Lambda}_{\mathrm{n}}}\mathbf{c}\label{eq:dispapprox final}
\end{align}
Taking the derivative of $\mathbf{u}(\xi)$ with respect to $\xi$ and substituting into~\Eref{eqn:sbfemstress} the stress field $\boldsymbol{\sigma}(\xi,\eta)$ can be expressed as:
\begin{align}
\boldsymbol{\sigma}(\xi,\eta)= & \boldsymbol{\Psi}_{\sigma}(\eta)\xi^{-\boldsymbol{\Lambda}_{\mathrm{n}}-\mathbf{I}}\mathbf{c}\label{eq:stress field complete}
\end{align}
where the stress mode $\boldsymbol{\Psi}_{\sigma}(\eta)$ is defined as: 
\begin{align}
\boldsymbol{\Psi}_{\sigma}(\eta)= & \mathbf{D}\left(-\mathbf{B}_{1}(\eta)\boldsymbol{\Phi}_{\mathrm{u}}\boldsymbol{\Lambda}_{\mathrm{n}}+\mathbf{B}_{2}(\eta)\boldsymbol{\Phi}_{\mathrm{u}}\right)\label{eq:stress mode}
\end{align}
The stiffness matrix of an element is obtained by first substituting \Eref{eq:int constants} into~\Eref{eqn:raddispfunc} at $\xi=1$. This results in:
\begin{align}
\mathbf{f}= & \boldsymbol{\Phi}_{\mathrm{q}}\boldsymbol{\Phi}_{\mathrm{u}}^{-1}\mathbf{\mathbf{u}_{\mathrm{b}}}\label{eq:equil poly}
\end{align}
From \Eref{eq:equil poly}, the stiffness matrix $\mathbf{K}$ can be identified to be given by the expression
\begin{equation}
\mathbf{K}=\boldsymbol{\Phi}_{\mathrm{q}}\boldsymbol{\Phi}_{\mathrm{u}}^{-1}\label{eqn:sbfemkmat-b}
\end{equation}

\begin{rmk}
The stiffness computed by employing the SBFEM is positive definite and symmetric. Hence, the stiffness matrix can be assembled in the conventional FEM approach. A simple Matlab \textsuperscript{\textregistered} function is given in~\cite{natarajanooi2014} to compute the stiffness matrix using the SBFEM.
\end{rmk}

\paragraph{Calculation of the stress intensity factors} A unique feature of the SBFEM is that stress singularities, if present,
are analytically represented in the radial displacement functions
$\mathbf{u}(\xi)$. When a crack is modelled by a polygon with its
scaling centre chosen at the crack tip in Figure (22), some of the
eigenvalues$\boldsymbol{\Lambda}_{\mathrm{n}}^{(\mathrm{s})}\subset\boldsymbol{\Lambda}_{\mathrm{n}}$
satisfy $-1<\boldsymbol{\Lambda}_{\mathrm{n}}^{(\mathrm{s})}<0$.
These eigenvalues lead to singular stresses at the crack tip. Using
$\boldsymbol{\Lambda}_{\mathrm{n}}^{(\mathrm{s})}$, the singular
stress field $\boldsymbol{\sigma}^{(\mathrm{s})}(\xi,\eta)$ can be
defined as {[}58{]}
\begin{align}
\boldsymbol{\sigma}^{(\mathrm{s})}(\xi,\eta)= & \boldsymbol{\Psi}_{\sigma}^{(\mathrm{s})}(\eta(\theta))\xi^{-\boldsymbol{\Lambda}_{\mathrm{n}}^{(\mathrm{s})}-\mathbf{I}}\mathbf{c}^{(\mathrm{s})}\label{eq:singstre}
\end{align}
where the singular stress mode $\boldsymbol{\Psi}_{\sigma}^{(\mathrm{s})}(\eta(\theta))=\left[\begin{array}[t]{ccc}
\boldsymbol{\Psi}_{\sigma_{xx}}^{(\mathrm{s})}(\eta(\theta))\quad & \boldsymbol{\Psi}_{\sigma_{yy}}^{(\mathrm{s})}(\eta(\theta))\quad & \boldsymbol{\Psi}_{\tau_{xy}}^{(\mathrm{s})}(\eta(\theta))\end{array}\right]^{\mathrm{T}}$ is
\begin{align}
\boldsymbol{\Psi}_{\sigma}^{(\mathrm{s})}(\eta(\theta))= & \mathbf{D}(-\mathbf{B}_{1}(\eta(\theta))\boldsymbol{\Phi}_{\mathrm{u}}^{(\mathrm{s})}\boldsymbol{\Lambda}_{\mathrm{n}}^{(\mathrm{s})}+\mathbf{B}_{2}(\eta(\theta))\boldsymbol{\Phi}_{\mathrm{u}}^{(\mathrm{s})})\label{eq:singstremode}
\end{align}
In \Eref{eq:singstremode} $\boldsymbol{\Phi}_{\mathrm{u}}^{(\mathrm{s})}\subset\boldsymbol{\Phi}_{\mathrm{u}}$ and $\mathbf{c}^{(\mathrm{s})}\subset\mathbf{c}$, contain the displacement modes and integration constants corresponding to $\boldsymbol{\Lambda}_{\mathrm{n}}^{(\mathrm{s})}$. It can be discerned from \Eref{eq:singstre} that $\boldsymbol{\Lambda}_{\mathrm{n}}^{(\mathrm{s})}$ leads to singular stresses at the crack tip. This enables the stress intensity factors to be computed directly from their definitions. The stress intensity factors for a crack that is aligned with the Cartesian coordinate axes shown in Figure (23) are defined as
\begin{align}
\left\{ \begin{array}{c}
K_{I}\\
K_{II}
\end{array}\right\} = & \stackrel{\mathrm{lim}}{r\rightarrow0}\left\{ \begin{array}{c}
\sqrt{2\pi r}\sigma_{yy}|_{\theta=0}\\
\sqrt{2\pi r}\tau_{xy}|_{\theta=0}
\end{array}\right\} \label{eq:SIF}
\end{align}
\begin{figure}[htpb]
\centering
\scalebox{0.7}{\input{./Figures/ScaledBFEMPoly.pstex_t}}
\caption{A cracked domain modelled by SBFEM and the definition of local coordinate system, where the `black' dots represent the nodes.}
\label{fig:crkPolySBFEM1}
\end{figure}
Substituting the stress components in \Eref{eq:singstre} at angle $\theta=0$ into Eq.~\eqref{eq:SIF} and using the relation $\xi=r/L_o$ at $\theta=0$, the stress intensity factors are
\begin{align}
\left\{ \begin{array}{c}
K_{I}\\
K_{II}
\end{array}\right\} = & \sqrt{2\pi L_{O}}\left\{ \begin{array}{c}
\boldsymbol{\Psi}_{\sigma_{yy}}^{(\mathrm{s})}(\eta(\theta=0))\mathbf{c}^{(\mathrm{s})}\\
\boldsymbol{\Psi}_{\tau_{xy}}^{(\mathrm{s})}(\eta(\theta=0))\mathbf{c}^{(\mathrm{s})}
\end{array}\right\} \label{eq:SIF2}
\end{align}

%% file: Figures/ScaledBFEMPoly.pstex_t
\begin{picture}(0,0)%
\includegraphics{./Figures/ScaledBFEMPoly.pstex}%
\end{picture}%
\setlength{\unitlength}{3947sp}%
\begingroup\makeatletter\ifx\SetFigFont\undefined%
\gdef\SetFigFont#1#2#3#4#5{%
  \reset@font\fontsize{#1}{#2pt}%
  \fontfamily{#3}\fontseries{#4}\fontshape{#5}%
  \selectfont}%
\fi\endgroup%
\begin{picture}(4685,4158)(3171,-5005)
\put(4906,-3781){\makebox(0,0)[lb]{\smash{{\SetFigFont{12}{14.4}{\familydefault}{\mddefault}{\updefault}$r_o$}}}}
\put(5116,-2191){\makebox(0,0)[lb]{\smash{{\SetFigFont{12}{14.4}{\familydefault}{\mddefault}{\updefault}$y$}}}}
\put(5746,-2386){\makebox(0,0)[lb]{\smash{{\SetFigFont{12}{14.4}{\familydefault}{\mddefault}{\updefault}$x$}}}}
\put(5071,-2776){\makebox(0,0)[lb]{\smash{{\SetFigFont{14}{16.8}{\familydefault}{\mddefault}{\updefault}$O$}}}}
\put(6316,-2671){\makebox(0,0)[lb]{\smash{{\SetFigFont{12}{14.4}{\familydefault}{\mddefault}{\updefault}$L_o$}}}}
\end{picture}%

%% file: Figures/dedgcrk.pstex_t
\begin{picture}(0,0)%
\includegraphics{./Figures/dedgcrk.pstex}%
\end{picture}%
\setlength{\unitlength}{3947sp}%
\begingroup\makeatletter\ifx\SetFigFont\undefined%
\gdef\SetFigFont#1#2#3#4#5{%
  \reset@font\fontsize{#1}{#2pt}%
  \fontfamily{#3}\fontseries{#4}\fontshape{#5}%
  \selectfont}%
\fi\endgroup%
\begin{picture}(3549,6654)(4151,-6703)
\put(7576,-3286){\rotatebox{90.0}{\makebox(0,0)[lb]{\smash{{\SetFigFont{14}{16.8}{\familydefault}{\mddefault}{\updefault}{\color[rgb]{0,0,0}$L=2$}%
}}}}}
\put(5176,-5086){\makebox(0,0)[lb]{\smash{{\SetFigFont{14}{16.8}{\familydefault}{\mddefault}{\updefault}{\color[rgb]{0,0,0}$H=1$}%
}}}}
\put(5281,-256){\makebox(0,0)[lb]{\smash{{\SetFigFont{14}{16.8}{\familydefault}{\mddefault}{\updefault}{\color[rgb]{0,0,0}$\sigma$}%
}}}}
\put(5371,-6616){\makebox(0,0)[lb]{\smash{{\SetFigFont{14}{16.8}{\familydefault}{\mddefault}{\updefault}{\color[rgb]{0,0,0}$\sigma$}%
}}}}
\put(6601,-3061){\makebox(0,0)[lb]{\smash{{\SetFigFont{14}{16.8}{\familydefault}{\mddefault}{\updefault}{\color[rgb]{0,0,0}$a$}%
}}}}
\put(4501,-3136){\makebox(0,0)[lb]{\smash{{\SetFigFont{14}{16.8}{\familydefault}{\mddefault}{\updefault}{\color[rgb]{0,0,0}$a$}%
}}}}
\end{picture}%

%% file: Figures/inclCrk.pstex_t
\begin{picture}(0,0)%
\includegraphics{./Figures/inclCrk.pstex}%
\end{picture}%
\setlength{\unitlength}{4144sp}%
\begingroup\makeatletter\ifx\SetFigFont\undefined%
\gdef\SetFigFont#1#2#3#4#5{%
  \reset@font\fontsize{#1}{#2pt}%
  \fontfamily{#3}\fontseries{#4}\fontshape{#5}%
  \selectfont}%
\fi\endgroup%
\begin{picture}(7500,6609)(1531,-6373)
\put(6241,-1846){\makebox(0,0)[lb]{\smash{{\SetFigFont{14}{16.8}{\familydefault}{\mddefault}{\updefault}$B$}}}}
\put(5086, 29){\makebox(0,0)[lb]{\smash{{\SetFigFont{14}{16.8}{\familydefault}{\mddefault}{\updefault}$\sigma_2$}}}}
\put(5281,-6286){\makebox(0,0)[lb]{\smash{{\SetFigFont{14}{16.8}{\familydefault}{\mddefault}{\updefault}$\sigma_2$}}}}
\put(5161,-4246){\makebox(0,0)[lb]{\smash{{\SetFigFont{14}{16.8}{\familydefault}{\mddefault}{\updefault}$2w$}}}}
\put(3841,-1501){\makebox(0,0)[lb]{\smash{{\SetFigFont{14}{16.8}{\familydefault}{\mddefault}{\updefault}$2w$}}}}
\put(1546,-2506){\makebox(0,0)[lb]{\smash{{\SetFigFont{14}{16.8}{\familydefault}{\mddefault}{\updefault}$\sigma_1$}}}}
\put(9016,-2551){\makebox(0,0)[lb]{\smash{{\SetFigFont{14}{16.8}{\familydefault}{\mddefault}{\updefault}$\sigma_1$}}}}
\put(5701,-2116){\makebox(0,0)[lb]{\smash{{\SetFigFont{14}{16.8}{\familydefault}{\mddefault}{\updefault}$\beta$}}}}
\put(5476,-3166){\makebox(0,0)[lb]{\smash{{\SetFigFont{14}{16.8}{\familydefault}{\mddefault}{\updefault}$2a$}}}}
\put(6796,-2566){\makebox(0,0)[lb]{\smash{{\SetFigFont{14}{16.8}{\familydefault}{\mddefault}{\updefault}$x_1$}}}}
\put(5626,-1396){\makebox(0,0)[lb]{\smash{{\SetFigFont{14}{16.8}{\familydefault}{\mddefault}{\updefault}$x_2$}}}}
\put(4651,-3451){\makebox(0,0)[lb]{\smash{{\SetFigFont{14}{16.8}{\familydefault}{\mddefault}{\updefault}$A$}}}}
\end{picture}%